\documentclass[12pt]{amsart}
\usepackage{euscript}
\usepackage{amssymb,amscd}
\newtheorem{theorem}{Theorem}[section]

\newtheorem{lemma}[theorem]{Lemma}
\newtheorem{proposition}[theorem]{Proposition}
\newtheorem{corollary}[theorem]{Corollary}

\newtheorem{remark}[theorem]{Remark}

\begin{document}

\title[Equilibrium Measures for Maps with Inducing Schemes]
{Equilibrium Measures for Maps with Inducing Schemes}

\author{Yakov Pesin}
\address{Department of Mathematics, Pennsylvania State
University, University Park, State College, PA 16802}
\email{pesin@math.psu.edu}
\author{Samuel Senti}
\address{Instituto de Matem\'atica, Universidade Federal do Rio de Janeiro, C.P. 68 530, CEP 21945-970, R.J., Brazil}
\email{senti@impa.br}

\date{\today}

\thanks{Y. Pesin was partially supported by the
National Science Foundation grant \#DMS-0503810, S. Senti was supported by the Swiss National Science Foundation.}
\subjclass{37D25, 37D35, 37E05, 37E10}

\begin{abstract}
We introduce a class of continuous maps $f$ of a compact topological space $I$ admitting inducing schemes and describe the tower constructions associated with them. We then establish a thermodynamic formalism, i.e., describe a class of real-valued potential functions $\varphi$ on $I$, which admit a unique equilibrium measure $\mu_\varphi$ minimizing the free energy for a certain class of invariant measures. We also describe ergodic properties of equilibrium measures including decay of correlation and the central limit theorem. Our results apply to certain maps of the interval with critical points and/or singularities (including some unimodal and multimodal maps) and to potential functions $\varphi_t=-t\log|df|$ with $t\in(t_0, t_1)$ for some $t_0<1<t_1$. In the particular case of $S$-unimodal maps we show that one can choose $t_0<0$ and that the class of measures under consideration consists of all invariant Borel probability measures.
\end{abstract}

\maketitle


\section{Introduction}


In this paper we develop a thermodynamic formalism for some classes of continuous maps of compact topological spaces. In the classical setting,  given a continuous map $f$ of a compact space $I$ and a continuous \emph{potential} function $\varphi$ on $I$, one studies the \emph{equilibrium} measures for $\varphi$, i.e., invariant Borel probability measures $\mu_\varphi$ on $I$ for which the supremum
\begin{equation}\label{sup}
\sup_{\mu\in\mathcal{M}(f,I)}\{h_\mu(f)+\int_I\,\varphi\,d\mu\}
\end{equation}
is attained, where $h_\mu(f)$ denotes the metric entropy and $\mathcal{M}(f,I)$ is the class of all $f$-invariant Borel probability measures on $I$. According to the classical variational principle (see for example, \cite{Ruelle2}) the above supremum is equal to the topological pressure $P(\varphi)$ of $\varphi$.

For a smooth one-dimensional map $f$ of a compact interval $I$ with critical points the ``natural'' class of potential functions includes functions which are not necessarily continuous such as, for instance, the function $\varphi(x)=-\log|df(x)|$ which is unbounded at critical points. Allowing noncontinuous potentials changes the setup, in particular, reducing the class of invariant measures under consideration. It also raises the question of adapting the notion of topological pressure to this new setup and establishing an appropriate version of the variational principle (we refer the reader to \cite{Mummert} for a discussion of these problems).

In this paper we develop a thermodynamic formalism for a class of maps admitting \emph{inducing schemes} satisfying some basic requirements. We establish ``verifiable'' conditions on potential functions, which guarantee the existence of a unique equilibrium measure for these potentials. We stress that one may have to restrict the supremum in \eqref{sup} to invariant measures satisfying some additional liftability requirements. Furthermore, the class of potential functions for which existence and uniqueness of equilibrium measures is guaranteed may depend on the choice of the inducing scheme.
Inducing schemes satisfying our requirements can be constructed for a broad class of one-dimensional maps, certain polynomial maps of the Riemann sphere and some multidimensional maps (see \cite{PSZ07}). We apply our results to study equilibrium measures for a broad class of one-dimensional maps (including $S$-unimodal maps) and for potential function $\varphi_t(x)=-t\log |df(x)|$ where $t$ runs some interval containing $[0,1]$.

In the first part, we describe an abstract \emph{inducing scheme} for a continuous map $f$ of a compact topological space of finite topological entropy. Such a scheme provides a symbolic representation of $f$, restricted to some invariant subset $X\subset I$, as a tower over $(W, F, \tau)$ where $F$ is the \emph{induced map} acting on the \emph{inducing domain} $W\subset I$ and $\tau$ is the \emph{inducing time}, which is a return time (not necessarily the first return time) to $W$. The level sets of the function $\tau$ are the basic elements of the inducing scheme. As the base $W$ of the tower can be a Cantor-like set, it can have a complicated topological structure.

An important feature of the inducing scheme is that basic elements form a countable generating Bernoulli partition for the induced map $F$ that is thus equivalent to the full shift on a countable set of states. Our results can be further generalized to towers for which the induced map $F$ is equivalent to a subshift of countable type provided it satisfies certain additional assumptions but we do not consider this case here.

The inducing procedures and the corresponding tower constructions for which the inducing time is the first return time to the base are classical objects in ergodic theory and were considered in works of Kakutani, Rokhlin, and others. Tower constructions for which the inducing time is not the first return time already appeared in works of J. Neveu \cite{Nev69} under the name of \emph{temps d'arret} and in the works of  Schweiger \cite{Schweiger75, Schweiger81} under the name \emph{jump transformation} (which are associated with some \emph{fibered systems}; see also the paper by Aaronson, Denker and Urba\'nski \cite{AaDU1} for some general results on ergodic properties of Markov fibered systems and jump transformations).

An $F$-invariant measure $\nu$ on $W$ with integrable inducing time (i.e., $\int_{W}\tau \,d\nu<\infty$) can be lifted to the tower thus producing an $f$-invariant measure $\mu={\mathcal L}(\nu)$ -- the \emph{lift} of $\nu$. Our thermodynamic formalism only allows $f$-invariant measures $\mu$ on $I$ that can be lifted. In particular, they should give positive weight to some invariant set $X\subseteq I$ (associated to the inducing scheme) which may be a proper subset of $I$. By Zweim\"uller \cite{Zweimuller7}, a measure $\mu$ on $X$ is liftable to the tower if it has integrable inducing time (i.e., $\int_X\tau\,d\mu<\infty$). The measure $\nu$ for which $\mu={\mathcal L}(\nu)$ is called the \emph{induced measure} for $\mu$ and is denoted by $i(\mu)$.

The liftability property is very important. In particular, for liftable measures one has Abramov's and Kac's formulas that connect respectively the entropy of the original map $f$ and the integral of the original potential $\varphi$ with the entropy of the induced map $F$ and the integral of the \emph{induced potential function} $\bar\varphi : W\to\mathbb{R}$ with respect to the induced measure. Whether a given invariant measure is liftable depends on the inducing scheme and there may exist non-liftable measures (see Section \ref{non-lift-example}). The \emph{liftability problem} is to construct, for a given map $f$, an ``optimal'' inducing scheme that captures \emph{all} invariant measures with positive weight to the base of the tower (i.e., every such measure is liftable). Such inducing schemes were studied in \cite{PSZ07}.

Our main result is that the lift of the equilibrium measure for the induced system is indeed the equilibrium measure for the original system. This is proven by studying the lift of a ``normalized'' potential cohomologous to $\varphi$. Also, we describe a condition on the potential function $\varphi$, which allows one to transfer results on ergodic properties of equilibrium measures for the induced system (including exponential decay of correlations and the Central Limit Theorem) to the original system. We stress again that the equilibrium measures we construct minimize the free energy $E_\mu=-h_\mu-\int_I\varphi\,d\mu$ only within the class of liftable measures, and we construct an example of an inducing scheme and a potential function $\varphi$ satisfying all our requirements, which possesses a unique non-liftable equilibrium measure (see Section \ref{non-lift-example}).

In the second part of the paper we apply our results to effect thermodynamic formalism for some one-dimensional maps. First, we present additional conditions on the inducing schemes namely bounded distortion, and a control of the size and number of the basic elements of a given inducing time (see Section~\ref{moreproperties}). These conditions are used in Section~\ref{jacobian} where we apply our results to one-dimensional maps and to the family of potential functions $\varphi_t(x)$ with $t$ in some interval $(t_0,t_1)$, $t_0< 1< t_1$. We establish existence and uniqueness of equilibrium measures (in the space of liftable measures). We also show how a sufficiently small exponential growth rate of the number of basic elements allows one to choose $t_0=-\infty$ and in particular,  to establish existence and uniqueness of the measure of maximal entropy (again within the class of liftable measures).

In this paper we are particularly interested in one principle example -- unimodal maps from a positive Lebesgue measure set of parameters in a transverse one-parameter family $f_a$ along with the potential function $\varphi_{t,a}(x)=-t\log|df_a(x)|$ with $t\in (t_0, t_1)$ for some $t_0<1<t_1$ (see Section~\ref{examples}). We show that the inducing scheme of~\cite{JCY3, Senti3} satisfies the slow growth rate condition on the number of basic partition elements of a given inducing time thus proving existence and uniqueness of equilibrium measures for $\varphi_{t,a}(x)$ for any $-\infty<t<t_1$ with $t_1=t_1(a)>1$. Applying results in \cite{Senti3} and \cite{Bru2}, we then solve the liftability problem in this case.

Our main result then claims that under the negative Schwarzian derivative assumption the inducing scheme of~\cite{JCY3, Senti3} is ``optimal'' in the sense that the supremum in \eqref{sup} can be taken over \emph{all} $f$-invariant Borel probability measures: For a transverse one-parameter family $f_a$ of unimodal maps of positive Lebesgue measure in the parameter space, there exists a unique equilibrium measure with respect to all (not only liftable) invariant measures associated to the potential function $\varphi_{t,a}(x)$ for any $t_0<t<t_1$ where $t_0=t_0(a)<0$ and $t_1=t_1(a)>1$. This extends the results of Bruin and Keller \cite{Bru-Kel} for the parameters under consideration. In particular, this also establishes the existence and uniqueness of the measure of maximal entropy by a different method than Hofbauer \cite{Hof1, Hof2}.

Finally, in Section~\ref{multimodalsection} we show that for potentials $\varphi_t(x)$ with $t$ close to~$1$ our results extend to some more general families of one-dimensional maps such as certain families of multimodal maps introduced by Bruin, Luzzatto and van Strien \cite{BLvS1} and cups maps as presented in \cite{DobbsPhD}.

Recently, Bruin and Todd \cite{BruTod07a} applied the results presented here (see also\cite{Pesin-Senti2}) to certain multimodal maps and prove the existence and uniqueness of equilibrium measure with respect to all invariant measures. They can deal with the liftability problem by building various inducing schemes and comparing the equilibrium measures associated to these schemes. The liftability problem for complex polynomials is also addressed in \cite{BruTod07}, and another class of potential functions is studied in \cite{BruTod07b}

By a recent result of Dobbs \cite{Dobbs07},  for the quadratic family, there exists a set of parameters $\mathcal{B}$ of positive Lebesgue measure such that for every $b\in\mathcal{B}$ one can find
$t_b\in (0,1)$ for which the phase transition occurs: the function
$\varphi_{t_b}$ possesses two equilibrium measures. We observe that the maps $f_b$ with $b\in\mathcal{B}$ are finitely (not infinitely) renormalizable while the unimodal maps for which our Theorem~\ref{mainunimodaltheorem} holds are non-renormalizable. At this point we pose the following problem:

\emph{Given a transverse family of $S$-unimodal maps, is there a set $\mathcal{A}$ of parameters of positive Lebesgue measure such that for every $a\in\mathcal{A}$ and every $t\in(-\infty,\infty)$ the function $\varphi_{t,a}$ possesses a unique equilibrium measure? Furthermore, is the pressure function $P(\varphi_{t,a})$ real analytic in~$t$?}

An affirmative solution of this problem would allow one, among other things, to further develop thermodynamic formalism for unimodal maps.

{\bf Structure of the paper.}
In Section~\ref{inducing}, we give a formal description of general inducing schemes. In Section \ref{symbolic} we state some results on existence and uniqueness of Gibbs (and equilibrium) measures for the one-sided Bernoulli shift (hence, for the induced map $F$) and for the induced potential; see Sarig \cite{Sarig1, Sarig2} and also Mauldin and Urba\'nski \cite{MauUrb1}, Yuri \cite{Yuri1} and Buzzi and Sarig \cite{BuzziSarig1}. In Section~\ref{thermodynamics} we introduce a set of conditions on the potential functions $\varphi$, which ensure that the corresponding induced potential functions $\overline{\varphi}$ possess unique equilibrium measures with respect to the induced system. These conditions are stated in terms of the inducing scheme and hence the class of potential functions to which our results apply depend on the choice of the inducing scheme. In section~\ref{moreproperties}, we provide some additional conditions on the inducing scheme, which then allow us to prove, in Section~\ref{jacobian}, that the potential functions $\varphi_t$ satisfy the conditions of Section~\ref{thermodynamics} for all $t_0<t<t_1$ with $t_0<0$ and $t_1>1$. In Section~\ref{examples}, we build an inducing scheme for a positive Lebesgue measure set of parameters in a one-parameter family of unimodal maps, which satisfy the conditions of Sections~\ref{inducing} and \ref{moreproperties}. We also address the liftability problem, proving that all measures of positive entropy, which give positive weight to the tower, are liftable. Moreover, we prove that measures of zero entropy and measures that are not supported on the tower cannot be equilibrium measures for $\varphi_t$ with $t_0<t<t_1$, thus proving existence and uniqueness of the equilibrium measure among all invariant measures. In Section~\ref{multimodalsection} we provide more examples, namely certain multimodal maps, cusp maps and one-dimensional complex polynomials.

{\bf Acknowledgments.} We would like to thank H. Bruin, J. Buzzi,
D. Dolgopyat,  F. Ledrappier, S. Luzzatto, M. Misiurewicz, O. Sarig,
M. Viana, M. Yuri and K. Zhang for valuable discussions and
comments. Finally, we thank the ETH, Z\"urich where part of this work
was conducted. Ya. Pesin wishes to thank the Research Institute for
Mathematical Science (RIMS), Kyoto and Erwin Schr\"odinger
International Institute for Mathematics (ESI), Vienna, where a part of
this work was carried out, for their hospitality.


\vskip.3in
\noindent
{\bf Part I: General Inducing Schemes.}


\section{Inducing Schemes and Their Properties}
\label{inducing}


Let $f:I\to I$ be a continuous map of a compact topological space $I$. Throughout this paper we shall always assume that the topological entropy $h(f)$ of $f$ is finite; in particular, the metric entropy $h_\mu(f)<\infty$ for any $f$-invariant Borel measure $\mu$. Let  $S$ be a countable collection of disjoint Borel subsets of $I$ called \emph{basic elements} and $\tau : S\to\mathbb{N}$ a positive integer-valued function. Define the \emph{inducing domain} by
$$
W:=\bigcup_{J\in S}\, J,
$$
the \emph{inducing time} $\tau : I\to\mathbb{N}$ by
$$
\tau(x):=
\begin{cases}
\tau(J), & x\in J,\ J\in S\\
0, & \mbox{otherwise}.
\end{cases}
$$
Let $\overline J$ denote the closure of the set $J$. We say that $f$ admits an \emph{inducing scheme} $\{S,\tau\}$ if the following conditions hold:
\begin{enumerate}
\item [(H1)]
for each $J\in S$ there exists a connected open set $U_J\supseteq J$ such that $f^{\tau(J)}|U_J$ is a homeomorphism onto its image and
$f^{\tau(J)}(J)=W$;
\item [(H2)] the partition $\mathcal{R}$ of $W$ induced by the sets $J\in S$ is Bernoulli generating in the following sense: for any countable collection of elements $\{J_k\}_{k\in\mathbb N}$, the intersection
$$
\overline{J_1}\cap\Bigl(\bigcap_{k\ge 2}
f^{-\tau(J_1)}\circ\cdots\circ f^{-\tau(J_{k-1})}(\overline{J_k})\Bigr)
$$
is not empty and consists of a single point, where $f^{-\tau(J)}$ denotes the inverse branch of $f^{\tau(J)}|J$ (here
$f^{-\tau(J)}(I)=\emptyset$ if $I\cap~f^{\tau(J)}(J)=~\emptyset$).
\end{enumerate}
Define the \emph{induced map} $F: W\to W$ by $F(x)=f^{\tau(x)}(x)$
and then set
\begin{equation}\label{tower}
X=\bigcup_{J\in S}\bigcup_{k=0}^{\tau(J)-1}f^k(J).
\end{equation}
The set $X$ is forward invariant under $f$. We also set
\begin{equation}\label{setw}
\mathcal W=\bigcup_{J\in S}\,\overline{J}.
\end{equation}
Conditions (H1) and (H2) allow one to obtain a symbolic representation of the induced map $F$ via the Bernoulli shift on a countable set of states. Consider the \emph{full shift of countable type} $(S^{\mathbb N},\sigma)$ where $S^{\mathbb N}$ is the space of one-sided infinite sequences with elements in $S$ and $\sigma$ is the (left) shift on $S^{\mathbb{N}}$, $(\sigma(a))_k:=a_{k+1}$ for $a=(a_k)_{k\ge 0}$. Define the \emph{coding map} $h\colon S^{\mathbb N}\to{\mathcal W}$ by $h((a_k)_{k\in\mathbb N}):=x$ where $x$ is such that $x\in\overline{J}_{a_0}$ and
\[
f^{\tau(J_{a_k})}\circ\cdots\circ f^{\tau(J_{a_0})}(x)\in
\overline{J}_{a_{k+1}} \quad\mbox{ for }\quad k\ge 0.
\]
\begin{proposition}\label{conjugation}
The map $h$ is well-defined, continuous and
$W\subseteq h(S^{\mathbb N})$. It is one-to-one on $h^{-1}(W)$ and is a conjugacy between $\sigma|h^{-1}(W)$ and $F|W$, i.e.,
$$
h\circ\sigma|h^{-1}(W)=F\circ h|h^{-1}(W).
$$
\end{proposition}
\begin{proof}
By (H2), given $a=(a_k)_{k\ge 0}$, there exists a unique point $x\in I$ such that $h(a)=x$. It follows that $h$ is well-defined. Moreover, given $x\in W$, there is a unique $a=(a_k)_{k\ge 0}$ such that
\[
f^{\tau(J_{a_k})}\circ\cdots\circ f^{\tau(J_{a_0})}(x)\in
J_{a_{k+1}} \quad\mbox{ for }\quad k\ge 0.
\]
It follows that $W\subseteq h(S^{\mathbb N})$ and that $h$ is one-to-one on $h^{-1}(W)$. Clearly, $\sigma|h^{-1}(W)$ and $F|W$ are conjugate via $h$. By (H2), for any $a=(a_k)_{k\ge 0}$ the sets $h([a_0,\dots, a_k])$ form a basis of the topology at $x=h(a)$. This implies that $h$ is continuous.
\end{proof}
Observe that the set $S^{\mathbb N}\setminus h^{-1}(W)$ contains no open subsets: indeed, by Conditions (H1) and (H2), the image of any cylinder $[a_1\dots a_n]$ under the coding map $h$ must contain points in $W$. This means that the set
$S^{\mathbb N}\setminus h^{-1}(W)$ is ``small'' in the topological sense but we also need it to be small in the measure-theoretical sense. More precisely, we require the following condition:
\begin{enumerate}
\item[(H3)] if $\mu$ is a shift invariant measure, which gives positive weight to any open set, then the set
$S^{\mathbb N}\setminus h^{-1}(W)$ has zero measure.
\end{enumerate}
This condition allows one to transfer shift invariant measures on $S^{\mathbb N}$, which give positive weight to open sets (in particular, Gibbs measures), to measures on $W$ invariant under the induced map.

Let $\mathcal{M}(F,W)$ be the set of $F$-invariant ergodic Borel
probability measures on $W$ and $\mathcal{M}(f,X)$ the set of
$f$-invariant ergodic Borel probability measures on $X$. For any
$\nu\in\mathcal{M}(F,W)$ set
$$
\displaystyle{
Q_\nu:=\sum_{J\in S}\tau(J)\nu(J). }
$$
If $Q_\nu<\infty$ we define the {\it lifted measure} ${\mathcal L}(\nu)$ on $I$
in the following way (see for instance \cite{dMvS1}): for any
measurable set $E\subseteq I$,
$$
{\mathcal L}(\nu)(E):=\frac1Q_\nu\sum_{J\in S}\,\sum_{k=0}^{\tau(J)-1}
\nu(f^{-k}(E)\cap J).
$$
The following result is immediate.
\begin{proposition}\label{liftedmeasure}
If $\nu\in\mathcal{M}(F,W)$ with $Q_\nu<\infty$ then
${\mathcal L}(\nu)\in\mathcal{M}(f,X)$ with ${\mathcal L}(\nu)(X)=1$ and
${\mathcal L}(\nu)|W\ll~\nu$.
\end{proposition}
We consider the class of measures
$$
\mathcal{M}_L(f,X):=\{\mu\in\mathcal{M}(f,X)\colon\text{ there exists } \nu\in\mathcal{M}(F,W),\ {\mathcal L}(\nu)=\mu\}.
$$
We call a measure $\mu\in~\mathcal{M}_L(f,X)$ \emph{liftable}.
It follows from Proposition \ref{liftedmeasure} that $\nu$ is uniquely
defined. We call $\nu$ the \emph{induced measure} for $\mu$ and we
write $\nu=:i(\mu)$. Observe that $Q_{i(\mu)}<\infty$ for any
$\mu\in\mathcal{M}_L(f,X)$.

Let $\varphi:I\to\mathbb{R}$ be a Borel function. In what follows we shall always assume that $\varphi$ is well-defined and is finite at every point $x\in{\mathcal W}$ (see \eqref{setw}) and we call $\varphi$ a \emph{potential}. We define the \emph{induced potential} $\overline\varphi: W \to\mathbb{R}$ by
\begin{equation}\label{ind_func1}
\overline\varphi(x):=\sum_{k=0}^{\tau(J)-1}\varphi(f^k(x))\quad
\text{ for }x\in J.
\end{equation}
We stress that the function $\varphi$ need not be continuous but in what follows we will require that the induced function $\bar\varphi$ is continuous in the topology of $W$.

Although the induced map $F$ may \emph{not} be the first return time map, Abramov's formula, connecting the entropies of $F$ and $f$, and Kac's formula, connecting the integrals of $\varphi$ and $\overline\varphi$, still hold (\cite[Proposition~2]{Nev69}, see also \cite{Zweimuller7} and, for related results, \cite{Keller3}).

\begin{theorem}[Abramov's and Kac's Formulae]\label{Kac}
Let $\nu\in\mathcal{M}(F,W)$ with $Q_\nu<\infty$. Then
$$
h_\nu(F)=Q_\nu\cdot h_{{\mathcal L}(\nu)}(f)<\infty.
$$
If $\int_W \overline\varphi\,d\nu$ is finite then
$$
-\infty<\int_W\overline\varphi\,d\nu=Q_\nu\cdot\int_X\varphi\,d{\mathcal L}(\nu)
<\infty.
$$
\end{theorem}
\begin{proof}
For the proof of Abramov's formula we refer to \cite{Zweimuller7} (recall that we require the topological entropy of $f$ to be finite). To prove Kac's formula, using the definition of ${\mathcal L}(\nu)$, we get
\begin{multline*}
\int_W\overline\varphi\, d\nu=\int_W\sum_{k=0}^{\tau(x)-1}\varphi(f^kx)\, d\nu(x)=
\sum_{J\in S}\sum_{k=0}^{\tau(J)-1}\int_J\varphi(f^kx)\,d\nu|J(x)\\
=\sum_{J\in S}\sum_{k=0}^{\tau(J)-1}\int_X\varphi(y)\,
d\nu(f^{-k}y\cap J)=Q_\nu\cdot\int_X\varphi\, d {\mathcal L}(\nu).
\end{multline*}
The desired result follows.
\end{proof}

We now prove that the space of liftable measures $\mathcal{M}_L(f,X)$ is non-empty. To this end we observe that $\mathcal{M}_L(f,X)\subseteq \mathcal{M}(f,X)$ and that $\mu(W)>0$ for any $\mu\in\mathcal{M}_L(f,X)$.

\begin{theorem}\label{entropyrelation}
Let $\mu\in\mathcal{M}(f,X)$ and $\tau\in L^1(X,\mu)$. Then $\mu\in\mathcal{M}_L(f,X)$ and
$$
h_{i(\mu)}(F)=Q_{i(\mu)}\cdot h_\mu(f)<\infty.
$$
In addition, if $\int_X\varphi\,d\mu$ is finite, then
$$
-\infty<\int_W\overline\varphi\,d i(\mu)=Q_{i(\mu)}\cdot\int_X\varphi\,d\mu <\infty.
$$
\end{theorem}
\begin{proof}
By \cite{Zweimuller7} (see also \cite{Bru2} for related results), there exists a measure $i(\mu)\in\mathcal{M}(F,W)$ such that $i(\mu)$ is absolutely continuous with respect to $\mu$, $Q_{i(\mu)}<\infty$, and ${\mathcal L}(i(\mu))=\mu$. Therefore, $\mu\in\mathcal{M}_L(f,X)$. To prove the other claims apply Theorem~\ref{Kac} to the measure $i(\mu)$. Since $h_\mu(f)<\infty$ (due to our assumption that the topological entropy of $f$ is finite) and ${\mathcal L}(i(\mu))=\mu$, we get
$$
h_{i(\mu)}(F)=Q_{i(\mu)}\cdot h_{{\mathcal L}(i(\mu))}(f)=Q_{i(\mu)}\cdot
h_\mu(f)<\infty.
$$
If $\int_X\varphi\,d\mu$ is finite, we get
$$
\int_W\overline\varphi\,di(\mu) =Q_{i(\mu)}\cdot\int_{X}\varphi\,d{\mathcal L}(i(\mu))
=Q_{i(\mu)}\cdot\int_X\varphi\,d\mu.
$$
This completes the proof of the theorem.
\end{proof}


\section{Thermodynamics Of Subshifts of Countable Type}\label{symbolic}


Consider the full shift $\sigma$ on $S^{\mathbb{N}}$ and let
$\Phi: S^{\mathbb{N}}\to\mathbb{R}$ be a continuous function (with respect to the discrete topology on $S^{\mathbb{N}}$). The
\emph{$n$-variation} $V_n(\Phi)$ is defined by
$$
V_n(\Phi):=\sup_{[a_0,\dots,a_{n-1}]}
\sup_{\omega,\omega'\in [a_0,\dots,a_{n-1}]}\{|\Phi(\omega)-\Phi(\omega')|\},
$$
where the \emph{cylinder set} $[a_0,\dots,a_{n-1}]$ consists of all
infinite sequences $\omega=(\omega_k)_{k\ge 0}$ with
$\omega_0=a_0,\ \omega_1=a_1,\dots,\ \omega_{n-1}=a_{n-1}$.

The \emph{Gurevich pressure} of $\Phi$ is defined by
\begin{equation}\label{gur}
P_G(\Phi):=\lim_{n\to\infty}\frac1n\log
\sum_{\sigma^n(\omega)=\omega}\exp(\Phi_n(\omega))1_{[a]}(\omega),
\end{equation}
where $a\in S$, $1_{[a]}$ is the characteristic function of the cylinder $[a]$ and
$$
\Phi_n(\omega):=\sum_{k=0}^{n-1}\Phi(\sigma^k(\omega)).
$$
It can be shown (see \cite{Sarig2}, \cite{Sarig4}) that if
$\sum\limits_{n\geq 2}V_n(\Phi)<\infty$ then the limit in~\eqref{gur} exists, does not depend on $a$, is never $-\infty$ and
$$
P_G(\Phi)=\lim\limits_{n\to\infty}
\frac{1}{n}\log\sum\limits_{\sigma^n(\omega)=\omega}\exp\Phi_n(\omega).
$$
A measure $\nu=\nu_{\Phi}$ is called a \emph{Gibbs measure} for $\Phi$ if there exist constants $C_1>0$ and $C_2>0$ such that for any cylinder set $[a_0,\dots , a_{n-1}]$ and any
$\omega\in [a_0,\dots ,a_{n-1}]$ we have
\begin{equation}\label{gibbs}
C_1\le\frac{\nu([a_0,\dots,a_{n-1}])}
{\exp\left(-nP_G(\Phi)+\Phi_{n}(\omega)\right)}\le C_2.
\end{equation}
Let $\mathcal{M}(\sigma)$ be the class of all $\sigma$-invariant ergodic Borel probability measures on $S^\mathbb{N}$.  A
$\sigma$-invariant measure $\nu_\Phi$ is said to be an
\emph{equilibrium measure} for $\Phi$ if
$-\int_{S^\mathbb{N}}\Phi d\nu_\Phi<\infty$ and
\begin{equation}
h_{\nu_\Phi}(\sigma)+\int\Phi d\nu_\Phi=
\sup_{\nu\in{\mathcal M}(\sigma):-\int_{S^\mathbb{N}}\,\Phi d\nu<\infty}\{h_\nu(\sigma)+\int\Phi d\nu\ \}.
\end{equation}
Note that unlike the classical case of subshifts of finite type the supremum above is taken only over the (restricted) class of measures $\nu$ for which $-\int_{S^\mathbb{N}}\Phi\,d\nu<\infty$.

A $\sigma$-invariant Gibbs measure $\nu$ for $\Phi$ is an equilibrium measure for~$\Phi$ provided $-\sum_{b\in S}\,\nu([b])\log \nu([b])<\infty$ (\cite{Bowen1}, see also \cite{Sarig1}). The following results establish the variational principle, and the existence and uniqueness of Gibbs and equilibrium measures for the full shift of countable type and for a certain class of potential functions. Various versions of these results were obtained by Mauldin and Urba\'nski \cite{MauUrb1}, by Sarig \cite{Sarig2}, \cite{Sarig3}, \cite{Sarig1} and by Yuri \cite{Yuri1} (see also \cite{Aaronson1} and \cite{BuzziSarig1}). In our presentation we follow \cite{Sarig2}, \cite{Sarig1}.

\begin{proposition}\label{sarig}
Assume that the potential $\Phi$ is continuous and that
$\sup_{\omega\in S^\mathbb{N}}\Phi<\infty$.
The following statements hold.
\begin{enumerate}
\item [1.]  If $\sum\limits_{n\geq 2}V_n(\Phi)<\infty$, then the variational principle for $\Phi$ holds:
$$
P_G(\Phi)=\sup_{\stackrel{\nu\in{\mathcal M}(\sigma)}{-\int_{S^\mathbb{N}}\,\Phi d\nu<\infty}}\{h_\nu(\sigma)+\int\Phi d\nu\ \}.
$$
\item [2.] If $\sum\limits_{n\ge 1}V_n(\Phi)<\infty$ and $P_G(\Phi)<\infty$, then there exists an ergodic $\sigma$-invariant Gibbs measure $\nu_{\Phi}$ for $\Phi$. If in addition, the entropy
$h_{\nu_\Phi}(\sigma)<\infty$, then $\nu_\Phi$ is the unique Gibbs and equilibrium measure.
\end{enumerate}
\end{proposition}
Observe that a Gibbs measure $\nu_\Phi$ is ergodic and positive on every non-empty open set.

In order to describe some ergodic properties of equilibrium measures
let us recall some definitions. A continuous transformation $T$ has
\emph{exponential decay of correlations} with respect to an invariant
Borel probability measure $\mu$ and a class $\mathcal{H}$ of functions if there exists $0<\theta<1$ such that, for any
$h_1, h_2\in\mathcal{H}$,
$$
\Big |\int h_1(T^n(x))h_2(x)\,d\mu -\int h_1(x)\,d\mu\int
h_2(x)\,d\mu\Big |\le K \theta^{n},
$$
for some $K=K(h_1,h_2)>0$.

The transformation $T$ satisfies the {\it central limit theorem} (CLT) for functions in $\mathcal{H}$ if for any $h\in\mathcal{H}$, which is not a coboundary (i.e., $h\ne g\circ T-g$ for any $g$), there exists $\gamma>0$ such that
$$
\mu\Bigl\{\frac{1}{\sqrt{n}}\sum_{i=0}^{n-1}(h(T^ix)-\int
h\,d\mu)<t\Bigr\}\rightarrow\frac{1}{\gamma\sqrt{2\pi}}\int_{-\infty}^t e^{-\tau^2/2\gamma^2}\,d\tau.
$$
The following statement describes ergodic properties of the equilibrium measure $\nu_\Phi$ and is a corollary of the well-known results by Ruelle \cite{Ruelle2} (see also \cite{Aaronson1}, \cite{Gordin} and \cite{Liverani1}). We say that the function $\Phi$ is \emph{locally
H\"older continuous} if there exist $A>0$ and $0<r<1$ such that for all $n\geq 1$,
\begin{equation}\label{clt}
V_n(\Phi)\leq A r^n.
\end{equation}

\begin{proposition}\label{sarig1}
Assume that $P_G(\Phi)<\infty$, $\sup_{\omega\in S^\mathbb{N}}\Phi<\infty$ and that $\Phi$ is locally H\"older continuous. If $h_{\nu_\Phi}(\sigma)<\infty$ then the measure $\nu_\Phi$ has exponential decay of correlations and satisfies the CLT with respect to the class of bounded H\"older continuous functions.
\end{proposition}


\section{Thermodynamics Associated with an Inducing Scheme}\label{thermodynamics}


\subsection{Classes of measures and potentials}
Let $f$ be a continuous map of a compact topological space $I$
admitting an inducing scheme $\{S,\tau\}$ satisfying conditions
(H1)--(H3) as described in Section~\ref{inducing}. Let also $\varphi:
X\to\mathbb{R}$ be a potential function, $\overline\varphi$ its
induced function, and $\mathcal{M}_L(f,X)$ the class of
liftable measures. We write
\begin{equation}\label{supremum}
P_L(\varphi):=\sup_{\mathcal{M}_L(f,X)}\{h_\mu(f)+\int_X\varphi\,d\mu\}
\end{equation}
and we call a measure $\mu_\varphi\in\mathcal{M}_L(f,X)$ an
\emph{equilibrium measure} for $\varphi$ (with respect to the class of measures $\mathcal{M}_L(f,X)$) if
$$
h_{\mu_\varphi}(f)+\int_X\varphi\,d\mu_\varphi=P_L(\varphi).
$$

Let us stress that our definition of equilibrium measures differs from
the classical one as we only allow \emph{liftable} measures, which give full weight to the \emph{noncompact} set $X$. Note that in general $P_L(\varphi)$ may not be finite and so we will need to impose conditions on the potential function in order to guarantee the finiteness of~$P_L(\varphi)$.

While dealing with the class of all $f$-invariant ergodic Borel
probability measures ${\mathcal M}(f,I)$,  depending on the potential
function $\varphi$, one may expect the equilibrium measure
$\mu_\varphi$ to be either non-liftable or to be supported outside of
the tower, i.e, $\mu_\varphi(X)=0$. In \cite{PesinZhang}, an example
of a one-dimensional map of a compact interval is given, which admits
an inducing scheme $\{S,\tau\}$ and a potential function $\varphi$
such that there exists a unique equilibrium measure $\mu_\varphi$ for
$\varphi$ (with respect to the class of measures ${\mathcal M}(f,I)$)
with $\mu_\varphi(X)=0$. The liftability problem is addressed in
\cite{PSZ07, PesinZhang1} where some characterizations of and criteria for liftability are obtained. Let us point out that non-liftable measures may exist and the liftability property of a given invariant measure depends on the inducing scheme. For certain interval maps, for instance, one can construct different inducing schemes over the same base such that a measure with positive weight to the base is liftable with respect to one of the schemes but not with respect to the other (see \cite{PesinZhang1}, also \cite{Bru2}). In Sections \ref{examples} and \ref{multimodalsection} we discuss liftability for unimodal and multimodal maps satisfying the Collet-Eckmann condition. In these particular cases we show that every measure in $\mathcal{M}(f,X)$ is liftable.

Two functions $\varphi$ and $\psi$ are said to be {\it cohomologous} if there exists a bounded function $h$ and a real number $C$ such that $\varphi-\psi=h\circ f-h+C$. An equilibrium measure for $\varphi$ is also an equilibrium measure for any $\psi$ cohomologous to $\varphi$. In particular, if $\varphi$ satisfies the conditions of Theorem~\ref{remark} below, then there exists a unique equilibrium measure for any $\psi$ cohomologous to $\varphi$ regardless of whether $\psi$ satisfies these conditions or not.

\subsection{Gibbs and equilibrium measures for the induced map}
In order to prove the existence of a unique equilibrium measure
$\nu_{\overline\varphi}$ for the induced map $F$ we impose some
conditions on the induced potential function $\overline\varphi$.

\begin{remark}\label{extension}
Note that in view of \eqref{ind_func1}, given $J\in S$, the function $\overline\varphi$ can be naturally extended to the closure $\overline J$. This means that the function $\Phi:=\overline\varphi\circ h$ is well-defined on $S^\mathbb{N}$ where $h$ is the coding map (see Proposition \ref{conjugation}).
\end{remark}
We call a measure $\nu_{\overline\varphi}$ on $W$ a \emph{Gibbs measure} for $\overline\varphi$ if the measure
$(h^{-1})_*\,\nu_{\overline\varphi}$ is a Gibbs measure for the function $\Phi$. We call $\nu_{\overline\varphi}$ an \emph{equilibrium measure} for $\overline\varphi$ if
$-\int_{W}\overline\varphi\,d\nu_{\overline\varphi}<\infty$ and
$$
h_{\nu_{\overline\varphi}}(F)+\int_W\overline\varphi\,d\nu_{\overline\varphi}=
\sup_{\stackrel{\nu\in\mathcal{M}(F,W)}{-\int_{W}\overline\varphi\,d\nu<\infty}}
\{h_\nu(F)+\int_W\overline\varphi\,d\nu\}.
$$
We say that the potential $\overline\varphi$
\begin{enumerate}
\item [(a)] has \emph{summable variations} if the function $\Phi$ has
summable variations, i.e.,
$$
\sum\limits_{n\ge 1}V_n(\overline\varphi\circ h)=
\sum\limits_{n\ge 1}V_n(\Phi)<\infty;
$$
\item [(b)] has \emph{finite Gurevich pressure} if
$P_G(\overline\varphi\circ h)=P_G(\Phi)<\infty$.
\end{enumerate}
Note that the image under the coding map $h$ of any periodic orbit for the shift $\sigma$ is a periodic orbit for the map $f$. Nevertheless, it may happen that the induced map $F$ possesses no periodic orbit. This is why from now on we assume that $F$ has at least one periodic orbit. In all interesting cases this requirement is satisfied.
\begin{theorem}\label{boundedenergy}
Assume that the function $\overline\varphi$ has summable variations and finite Gurevich pressure. Then
$$
-\infty<P_L(\varphi)<\infty.
$$
\end{theorem}

\begin{proof}
By Proposition~\ref{conjugation}, there is a periodic orbit for $F$ in the set $W$. For the Dirac measure on that orbit we have that $\int_X\varphi\,d\mu>-\infty$. Since $0\le h_\mu(f)$, we conclude that $P_L(\varphi)>-\infty$.

For every $\mu\in\mathcal{M}_L(f,X)$ there exists a measure $i(\mu)\in\mathcal{M}(F,W)$ with $Q_{i(\mu)}<\infty$ and by Theorem~\ref{Kac},
\begin{equation}\label{boundedentropy}
0\le~h_{i(\mu)}(F)=Q_{i(\mu)}\cdot h_\mu(f)<\infty.
\end{equation}
Take $\mu\in\mathcal{M}_L(f,X)$ such that
$\int_W\overline\varphi\,di(\mu)>-\infty$. Since
$\overline\varphi$ has summable variations and finite Gurevich
pressure (and is thus bounded from above), Proposition \ref{sarig} and
\eqref{boundedentropy} imply that
$-\infty<\int_W\overline\varphi\,di(\mu)<\infty$ and, by Theorem \ref{Kac},
$$
-\infty<\int_W\overline\varphi\,di(\mu)=Q_{i(\mu)}\cdot
\int_X\varphi\,d\mu<\infty.
$$
If $P_L(\varphi)$ is non-positive the upper bound is immediate.
If $P_L(\varphi)$ is positive, using the fact that $1\le Q_{i(\mu)}<\infty$, we get
\begin{multline*}
P_L(\varphi)=
\sup_{\stackrel{\mu\in{\mathcal M}_L(f,X)}{\int_W\overline\varphi\,di(\mu)>-\infty}}
\left(\frac{h_{i(\mu)}(F)+\int_W\overline\varphi\,di(\mu)}{Q_{i(\mu)}}\right)\\
\le\,
\sup_{\stackrel{\nu\in{\mathcal M}(F,W)}{-\int_{W}\overline\varphi\,d\nu<\infty}}
\left(h_\nu(F)+\int_W\overline\varphi\,d\nu\right)<\infty,
\end{multline*}
where the first equality follows from the fact that $P_L(\varphi)$ cannot
be achieved by a measure with
$\int_W\overline\varphi\,di(\mu)=-\infty$. Indeed otherwise,
\begin{multline*}
\int_X\varphi(x)\,d{\mathcal L}(i(\mu))(x)=
\int_X\frac{1}{Q_{i(\mu)}}\varphi(x)\sum_{J\in S}\sum_{k=0}^{\tau(J)-1}\,di(\mu)(f^{-k}(x)\cap J)\\
=\frac{1}{Q_{i(\mu)}}\int_W\sum_{J\in S}\sum_{k=0}^{\tau(J)-1}\varphi(f^k(y)) di(\mu)(y\cap J)\\
=\frac{1}{Q_{i(\mu)}}\int_W\overline\varphi(y)\,di(\mu)=-\infty
\end{multline*}
would imply $P_L(\varphi)=-\infty$ contradicting the lower bound
established above.
\end{proof}

In order to show that equilibrium measures for the induced system lift to equilibrium measures for the original system, it is useful to work with a potential function which is cohomologous to the original potential function $\varphi$: when $P_L(\varphi)$ is finite we denote the induced function for $\varphi-P_L(\varphi)$ by $\varphi^+:=\overline{\varphi-P_L(\varphi)}=\overline\varphi-~P_L(\varphi)\tau$.
Given $J\in S$, this function can be naturally extended to the closure $\overline J$ and hence the function $\Phi^+:=\varphi^+\circ h$ is well-defined on $S^\mathbb{N}$ where $h$ is the coding map (see Proposition \ref{conjugation}). The following statement establishes the existence and uniqueness of equilibrium measures for $\varphi^+$ for the induced map $F$.

\begin{theorem}\label{arprinc3}
Assume that the induced function $\overline\varphi$ on $W$ has summable variations and finite Gurevich pressure. Also assume that  the function $\varphi^+$ has finite Gurevich pressure and that
\begin{equation}\label{bound1}
\sup_{J\in S}\sup_{x\in{\overline J}}\,\varphi^+(x)<\infty.
\end{equation}
Then the following statements hold:
\begin{enumerate}
\item [1.]  there exists an $F$-invariant ergodic Gibbs measure $\nu_{\varphi^+}$ on $W$, which is unique when
$h_{\nu_{\varphi^+}}(F)<\infty$;
\item [2.] if $Q_{\nu_{\varphi^+}}<\infty$ then
$\nu_{\varphi^+}$ is the unique equilibrium measure among the measures $\nu\in\mathcal{M}(F,W)$ satisfying
$\int_W\overline\varphi\,d\nu>-\infty$.
\end{enumerate}
\end{theorem}
\begin{proof}
Since $\overline\varphi$ has summable variations, it is continuous on $W$. Note that the inducing time $\tau$ is constant on
elements $J\in S$. It follows that the function $\varphi^+$ is also
continuous on $W$ and has summable variations. In view of \eqref{bound1}, we can apply Proposition~\ref{sarig} proving the existence of a $\sigma$-invariant ergodic Gibbs measure for $\Phi^+$. As a Gibbs measure must give positive weight to cylinders, it cannot be supported on $S^{\mathbb{N}}\setminus h^{-1}(W)$ due to Condition (H3) and the first statement follows.

For an $f$-invariant Borel probability measure $\mu$, we have
$0\le h_\mu(f)<\infty$. Theorem~\ref{Kac} and the assumption
$Q_{\nu_{\varphi^+}}<\infty$ imply $h_{\nu_{\varphi^+}}(F)<\infty$.
The second statement then follows from Proposition~\ref{sarig}.
\end{proof}

\subsection{Lifting Gibbs measures} We now describe a condition on the induced function $\overline\varphi$, which will help us prove that the \emph{natural candidate} -- the lifted measure
$\mu_\varphi:={\mathcal L}(\nu_{\varphi^+})$ where the measure $\nu_{\varphi^+}$ is constructed in Theorem~\ref{arprinc3} -- is indeed an equilibrium measure for $\varphi$.

We say that the induced function $\overline{\varphi}$ is
\emph{positive recurrent} if there exists $\varepsilon_0>0$ such that
$$
\varphi^+_{\varepsilon_0}:=\overline{\varphi-P_L(\varphi)+\varepsilon_0}
=\varphi^++\varepsilon_0\tau
$$
has finite Gurevich pressure. It follows that for any $0\le\varepsilon\le\varepsilon_0$ the function $\varphi^+_\varepsilon:=\overline{\varphi-P_L(\varphi)+\varepsilon}=\varphi^++\varepsilon\tau$
also has finite Gurevich pressure.
\begin{theorem}\label{arprinc1}
Assume that the induced function $\overline\varphi$ on $W$ has summable variations, finite Gurevich pressure and is positive recurrent. Also assume that the function $\varphi^+$ satisfies \eqref{bound1} and $Q_{\nu_{\varphi^+}}<\infty$ for the equilibrium measure $\nu_{\varphi^+}$ of theorem~\ref{arprinc3}. Then the measure $\mu_\varphi={\mathcal L}(\nu_{\varphi^+})$ is the unique equilibrium measure for $\varphi$ with respect to the class of liftable measures $\mathcal{M}_L(f,X)$.
\end{theorem}
\begin{proof}
Since $\overline\varphi$ is positive recurrent, the function
$\varphi^+$ has finite Gurevich pressure and all requirements of  Theorem~\ref{arprinc3} hold. By this theorem, the measure
$\mu_\varphi$ is well defined and belongs to $\mathcal{M}_L(f,X)$. We show that $P_G(\Phi^+)=0$ and that $\mu_\varphi$ is the unique equilibrium measure (with respect to the class of measures
$\mathcal{M}_L(f,X)$). As
$h_{\mu_{\varphi}}(f)+\int_X(\varphi-P_L(\varphi))\,d\mu_{\varphi}\le 0$ and $Q_{\nu_{\varphi^+}}\in[1,\infty)$, Proposition~\ref{sarig} and
Theorem~\ref{Kac} imply
\begin{multline}\label{Pis0}
P_G(\Phi^+)=h_{\nu_{\varphi^+}}(F)
+\int_W\varphi^+\,d\nu_{\varphi^+}\\
=Q_{\nu_{\varphi^+}}\cdot\bigl(h_{\mu_{\varphi}}(f)+\int_X(\varphi-P_L(\varphi))\,d\mu_{\varphi}\bigr)\le 0.
\end{multline}
On the other hand, for every $\varepsilon>0$ there is
$\mu\in\mathcal{M}_L(f,X)$ such that
$$
h_\mu(f)+\int_X\varphi\,d\mu\ge P_L(\varphi)-\varepsilon.
$$
Since $Q_{i(\mu)}$ is strictly positive for all $\mu$, Theorem~\ref{Kac} gives
\begin{multline*}
P_G(\Phi^+_\varepsilon)\ge h_{i(\mu)}(F)+\int_W\varphi^+_\varepsilon\,d i(\mu)\\
=Q_{i(\mu)}\cdot
\bigl( h_\mu(f)+\int_X(\varphi-P_L(\varphi)+\varepsilon)\,d\mu\bigr)\ge 0.
\end{multline*}
By \eqref{gur} and positive recurrence $P_G(\Phi^+_\varepsilon)$ is continuous in $\varepsilon$ for $0\le \varepsilon\le \varepsilon_0$. We conclude that $P_G(\Phi^+)\ge 0$, hence, \eqref{Pis0} becomes
\[
0=P_G(\Phi^+)=Q_{\nu_{\varphi^+}}\cdot
\bigl(h_{\mu_{\varphi}}(f)+\int_X(\varphi-P_L(\varphi))\,d\mu_{\varphi}\bigr).
\]
As $Q_{\nu_{\varphi^+}}\in[1,\infty)$, the measure $\mu_\varphi$ is an equilibrium measure for $\varphi$ (for the class of measures $\mathcal{M}_L(f,X)$). Unicity (over this class) follows from the unicity of $\nu_{\varphi^+}$.
\end{proof}

\subsection{Conditions on potential functions}\label{cond_pot_fun}
Verifying the hypotheses of Theorems~\ref{arprinc1} and \ref{ergodic} may be intricate. Additional conditions on the induced potential $\overline\varphi$ can help us check them.

Given a cylinder $[a_0,\dots,a_{n-1}]$, we set
$$
\begin{aligned}
J_{[a_0,\ldots,a_{n-1}]}:=&h([a_0,\dots,a_{n-1}])\\
=&\overline{J_{a_0}}\cap\Bigl(\bigcap_{k=2}^{n-1}
f^{-\tau(J_{a_0})}\circ\cdots\circ f^{-\tau(J_{a_{n-2}})} (\overline{J_{a_{n-1}}})\Bigr)
\end{aligned}
$$
(see Proposition~\ref{conjugation} for the definition of the conjugacy $h$). The \emph{$n$-variation} of $\overline\varphi$ is defined by
$$
V_n(\overline\varphi):=\sup_{[a_0,\dots,a_{n-1}]}
\sup_{x,x'\in J_{[a_0,\ldots,a_{n-1}]}}\{|\overline\varphi(x)-\overline\varphi(x')|\}.
$$
We assume the following conditions on the potential function $\varphi$:

\begin{enumerate}
\item [(P1)] $\overline\varphi$ is locally H\"older continuous (see \eqref{clt}): there exist $A>0$ and $0<r<1$ such that for all $n\geq 1$,
$$
V_n(\overline\varphi)\leq A r^n;
$$
\item [(P2)]
$$
\sum_{J\in S}\,\sup_{x\in J}\,\exp\,\overline\varphi(x)<\infty;
$$
\item [(P3)] there exists $\varepsilon_0>0$ such that
$$
\sum_{J\in S}\,\tau(J)\sup_{x\in J}\,\exp\,(\varphi^+(x)+\varepsilon_0 \tau(x))<\infty;
$$
\end{enumerate}

Let $\varphi$ be a bounded Borel function on $I$, which is H\"older continuous on the closure $\overline J$ of each $J\in S$.
Then $\varphi$ has bounded variation and there exists $C\ge 0$ such that the function $\varphi-c$ satisfies Condition (P2) for every $c\ge C$.

\begin{theorem}\label{remark}
Let $f$ be a continuous map of a compact topological  space. Assume that the topological entropy $h(f)<\infty$ and that $f$ admits an inducing scheme $\{S,\tau\}$ satisfying Conditions (H1)--(H3). Let $\varphi$ be a potential function satisfying Conditions (P1)--(P3). Then there exists a unique equilibrium measure $\mu_\varphi$ for $\varphi$ (with respect to the class of measures $\mathcal{M}_L(f,X)$).
\end{theorem}

\begin{proof} The proof will follow from Theorem \ref{arprinc1} if we prove that the induced potential $\overline\varphi$ satisfies its assumptions. By Condition (P1), the induced potential function
$\overline\varphi$ is continuous on $W$ and has summable variations.
Proposition~\ref{conjugation} implies that given any cylinder
$[a_0,\dots,a_{n-1}]$, there exists a unique $x$ in
$J_{[a_0,\dots,a_{n-1}]}$ with $F^n(x)=x$. Therefore, Condition (P2) implies
\begin{multline*}
\lim_{n\to\infty}\frac1n\log
\sum_{\stackrel{F^n(x)=x}{x\in J_{a_0}}}
\exp\left(\sum_{i=0}^{n-1}\overline\varphi(F^i(x))\right)\\
\le\lim_{n\to\infty}\frac1n\log\left(\sum_{J\in S}\,
\sup_{x\in J}\,\exp\,\overline\varphi(x)\right)^n<\infty,
\end{multline*}
thus proving that $\overline\varphi$ has finite Gurevich pressure. Positive recurrence follows from (P3) in the same way. Condition (P1) also implies that the induced function $\overline\varphi$ satisfies \eqref{bound1}. Together with Theorem~\ref{boundedenergy} this implies the finiteness of $P_L(\varphi)$, and so Conditions (P1) and (P3) (with $\varepsilon=0$) imply that the induced potential $\varphi^+$ corresponding to the ``normalized'' potential
$\varphi-P_L(\varphi)$ has summable variations and finite Gurevich pressure. By Theorem~\ref{arprinc3}, there exists a Gibbs measure
$\nu_{\varphi^+}$ for $\varphi^+$ on $W$. By \eqref{gibbs}, there
exist $C_1, C_2>0$ such that for every $J\in S$ and $x\in J$,
\begin{equation}\label{int-eq1}
C_1\le\frac{\nu_{\varphi^+}(J)}{\exp(-P+\varphi^+(x))}\le C_2,
\end{equation}
where $P=P_G(\Phi^+)$ is the Gurevich pressure of $\Phi^+$.
Summing \eqref{int-eq1} over all $J\in S$ and using Condition (P3) we get
$$
Q_{\nu_{\varphi^+}}=\sum_{J\in S}\tau(J)\nu_{\varphi^+}(J)\le
\frac{C_2}{e^P}\sum_{J\in S}\,\tau(J)\,\sup_{x\in J}\,
\exp\,(\varphi^+(x))<\infty.
$$
By Theorem~\ref{arprinc3}, $\nu_{\varphi^+}\in\mathcal{M}(F,W)$ is the unique equilibrium measure for $\varphi^+$ and, by Theorem~\ref{arprinc1}, ${\mathcal L}(\nu_{\varphi^+})\in\mathcal{M}_L(f,X)$ is the unique equilibrium measure (with respect to the class of measures $\mathcal{M}_L(f,X)$).
\end{proof}

\subsection{Ergodic properties}
To describe some ergodic properties of equilibrium measures we introduce another condition. Let $\varphi$ be a potential function. Consider the function $\varphi^+=\overline{\varphi-P_L(\varphi)}$ and let $\nu_{\varphi^+}$ be its equilibrium measure. We say that it has \emph{exponential tail} if there exist $K>0$ and $0<\theta<1$ such that for all $n>0$,
\begin{enumerate}
\item[(P4)] $\nu_{\varphi^+}(\{x\in W: \tau(x)\ge n\})\le K\theta^n.$
\end{enumerate}

\begin{theorem}\label{ergodic}
Assume that the induced function $\overline\varphi$ on $W$ is locally H\"older continuous, positive recurrent and has finite Gurevich pressure. Also assume that the function $\varphi^+$ satisfies Condition \eqref{bound1}. If $\nu_{\varphi^+}$ has exponential tail then there exists a unique equilibrium measure $\mu_\varphi$ (with respect to the class of measures $\mathcal{M}_L(f,X)$). It is ergodic, has exponential decay of correlations and satisfies the central limit theorem with respect to the class of functions whose induced functions on $W$ are bounded H\"older continuous.
\end{theorem}
\begin{proof}
If $\overline\varphi$ is locally H\"older continuous then it has
summable variations. Theorem~\ref{arprinc3} then implies the existence of a Gibbs measure $\nu_{\varphi^+}$. Since
$\nu_{\varphi^+}$ has exponential tail, we obtain
$$
Q_{\nu_{\varphi^+}}= \sum_{J\in S}\tau(J)\nu_{\varphi^+}(J)\le
\sum_{\ell=1}^\infty\ell\,\sum_{\stackrel{J\in S}{\tau(J)=\ell}}
\,\nu_{\varphi^+}(J)\le
K\,\sum_{\ell=1}^\infty\ell\,\theta^{\ell}<\infty.
$$
Since $\overline\varphi$ is positive recurrent, by
Theorem~\ref{arprinc1}, the measure $\mu_\varphi={\mathcal L}(\nu_{\varphi^+})$ is the unique equilibrium measure for $\varphi$. The desired result then follows from Theorem~\ref{sarig1} and Theorems 2 and 3 of Young in~\cite{LSY1}.
\end{proof}

\subsection{Non-liftable equilibrium measures}\label{non-lift-example}
We present an example of an inducing scheme $\{S,\tau\}$ for an interval map $f$ and a potential function $\varphi$ such that:
(1) $\varphi$ satisfies Conditions (P1)--(P3);
(2) $\varphi$ admits a unique equilibrium measure $\mu_\varphi$ within the class of all invariant measures, which gives positive weight to the base of the tower;
(3) $\mu_\varphi$ is not liftable.
Of course, by Theorem \ref{remark}, there exists another invariant measure, which is a unique equilibrium measure within the class of liftable measures.

Consider the map $f=2x \pmod 1$ of the unit interval $I$. The Lebesgue measure $\text{Leb}$ is the unique equilibrium measure of maximal entropy, i.e., the unique equilibrium measure for the potential function $\varphi=$constant.

Set $I^{(1)}=[0,\frac12 ]$, $I^{(2)}=(\frac12, 1]$ and consider the inducing scheme $\{S',\tau'\}$ where $S'$ is the countable collection of intervals $I_n$ such that $I_0=I^{(2)}$ and $I_n=f^{-1}(I_{n-1})\cap I^{(1)}$ for $n\ge 1$, and $\tau'(I_n)=n$. It is easy to see that this inducing scheme satisfies Conditions (H1)--(H3) and that the function $\varphi=-2$ satisfies Conditions (P1)--(P3) (with respect to the scheme $\{S',\tau'\}$). The corresponding equilibrium measure $\mu_\varphi=Leb$. In fact, every measure $\mu\in\mathcal{M}(f,X)$ is liftable to $\{S',\tau'\}$.

Now subdivide each interval $I_n$ into $2^{2^n}$ intervals of equal length and call them $I_n^j$. Consider the inducing scheme $\{S,\tau\}$ where $S$ consists of intervals $I_n^j$, $j=1,\dots, 2^{2^n}$, $n\ge 1$ and $\tau(I_n^j)=2^n+n$. It is shown in \cite{PesinZhang1} that $Leb$ is not liftable to $\{S,\tau\}$, however, it is easy to check that the function $\varphi=-2$ satisfies Conditions (P1)--(P3) (with respect to the inducing scheme $\{S,\tau\}$). By Theorem \ref{remark}, the function $\varphi$ possesses a unique equilibrium measure (within the class of liftable measures) $\mu_\varphi$, which is singular with respect to $Leb$.

In \cite{PesinZhang} the authors also provide examples of inducing schemes such that the supremum $P_L(\varphi)$ of \eqref{supremum} is strictly less than the supremum in \eqref{sup}.

The liftability problem for general piecewise invertible maps is addressed in detail in \cite{PSZ07}, see also \cite{BruTod07a} where the problem of comparing equilibrium measures obtained by different inducing schemes is addresses for certain multimodal maps and for the potentials $-t\log|df(x)|$ with $t$ close to $1$, see also \cite{BruTod07, BruTod07b}.

\vskip.3in\noindent
{\bf Part II: Applications to Interval Maps.}

\section{Inducing Schemes with Exponential Tail and Bounded Distortion}\label{moreproperties}


In this section we apply the above results to effect the thermodynamic formalism for $C^1$ maps $f$ of a compact interval $I$ that admit inducing schemes $\{S, \tau\}$. We shall study equilibrium measures corresponding to the special family of potential functions
$\varphi_t(x)=-t\log|df(x)|$ where $t$ runs in some interval of
$\mathbb{R}$. We shall show that $\varphi_t(x)$ satisfies Conditions
(P1)--(P4) of Part I for $t$ in some interval $(t_0,t_1)$ provided that the inducing scheme satisfies some additional properties, namely an exponential bound on the ``size'' of the partition elements with large inducing time, bounded distortion and a bound on the cardinality of partition elements with given inducing time. We also present some examples of systems, which admit such inducing schemes.

Denote the Lebesgue measure of the set $J\in S$ by $\text{Leb}(J)$. We assume that the inducing scheme satisfy the following additional conditions
\begin{enumerate}
\item [(H4)] \emph{exponential tail}: We have $\text{Leb}(W)>0$ and there are constants $c_1>0$ and $\lambda_1>1$ such that for all $n\ge 0$,
$$
\sum_{J\in S\colon\tau(J)\ge n}\text{Leb}(J)\le c_1^{-1}\lambda_1^{-n};
$$
\item [(H5)] \emph{bounded distortion}: there are constants $c_2>0$ and $\lambda_2>1$ such that for all $n\ge 0$, each cylinder $[a_0,\ldots,a_{n-1}]$, any two points $x, y\in J_{[a_0,\ldots,a_{n-1}]}$ (see Section~\ref{cond_pot_fun} for the definition of the set), and each $0\le i\le n-1$, we have
$$
\left|\frac{dF(F^i(x))}{dF(F^i(y))}-1\right|\le c_2\lambda_2^{-n}.
$$
\end{enumerate}
Conditions (H4) and (H5) imply the following.
\begin{corollary}\label{distbounds}
There are positive constants $c_3, c_4$ and $\lambda_3>1$ such that for every $J\in S$ and $x\in J$,
$$
c_1\lambda_1^{\tau(J)}\le\text{Leb}(J)^{-1}\le c_3|dF(x)|\le c_4\lambda_3^{\tau(J)}.
$$
\end{corollary}
\begin{proof}
The first inequality follows from (H4). Since $W=F(J)$ for any
$J\in S$, the other inequalities follow from Conditions (H4) and (H5) and the fact that the derivative is bounded from above on a compact interval $I$.
\end{proof}
\begin{remark}\label{refinement}
Without loss of generality one can assume that $c_1=1$. Indeed, partition elements of lower order can be refined and the constant $\lambda_1$ can be adjusted for this purpose. Obviously, one can also choose $\lambda_3$ such that $c_4=1$.
\end{remark}

\begin{theorem}\label{integrableN}
Assume that $f$ admits an inducing scheme $\{S, \tau\}$ satisfying Conditions~(H1)--(H5). Then for any measure $\mu\in\mathcal{M}_L(f,X)$,
\[
\log\lambda_1\le\int_X\,\log|df|\,d\mu\le\log\lambda_3.
\]
\end{theorem}
\begin{proof}
By Corollary~\ref{distbounds}, we have that for every $J\in S$ and any $x\in J$,
\begin{equation}\label{int_estimate1}
\tau(J)\log\lambda_1\le\log|dF(x)|\le\tau(J)\log\lambda_3.
\end{equation}
For any $\mu\in\mathcal{M}_L(f,X)$ integrating~\eqref{int_estimate1}
against $i(\mu)$ over $J$ and summing over all $J\in S$ yields
\[
Q_{i(\mu)}\log\lambda_1
\le\int_W\,\log|dF(x)|\,di(\mu)\le Q_{i(\mu)}\log\lambda_3.
\]
By Theorem~\ref{entropyrelation}, we have
\[
\int_W\,\log|dF(x)|\,di(\mu)=Q_{i(\mu)}\int_X\,\log|df(x)|\,d\mu,
\]
and the statement follows since $Q_{i(\mu)}$ is positive.
\end{proof}

As an immediate corollary of this result we obtain the following statement.

\begin{corollary}\label{integrableN1}
Assume that $f$ admits an inducing scheme $\{S, \tau\}$ satisfying
Conditions~(H1)--(H5). Then for any ergodic measure $\mu\in\mathcal{M}_L(f,X)$ the Lyapunov exponent $\lambda(\mu)$ of $\mu$ is strictly positive. Moreover,
$\log\lambda_1\le\lambda(\mu)\le\log\lambda_3$.
\end{corollary}

\begin{proof}
It suffices to notice that $\lambda(\mu)=\int_X\,\log|df|\,d\mu$ and
use Theorem~\ref{integrableN}.
\end{proof}

Denote by $S(n):=\text{Card}\{J\in S\ |\ \tau(J)=n\}$. Conditions~(H4) and (H5) imply that
\begin{equation}\label{sng}
S(n)\le c_6\gamma^n
\end{equation}
for some $1\le\gamma\le \frac{\lambda_3}{\lambda_1}$ and
$c_6=c_6(\gamma)>0$. For our main results we need a better control of the growth rate of $S(n)$, which is given by the following condition
\begin{enumerate}
\item[(H6)] \emph{Subexponential growth of basic elements}: for every $\gamma>1$ there exists $d>0$ such that $S(n)\le d\gamma^n$ for every $n\ge 1$.
\end{enumerate}


\section{Equilibrium Measures For Potentials $-t\log|df(x)|$}
\label{jacobian}


We now apply the results of the previous sections to the family of
potential functions $\varphi_t(x)=-t\log |df(x)|$, $x\in I$ for
$t\in\mathbb{R}$. The corresponding induced potential is
$$
\overline\varphi_t(x)=\sum_{k=0}^{\tau(x)-1}-t\log|df(f^k(x))|=
-t\log|dF(x)|.
$$
Given $c\in\mathbb{R}$, we also consider the \emph{shifted} potential
$\xi_{c,t}:=\varphi_t+c$ and its \emph{induced potential}
$$
\overline{\xi}_{c,t}(x):=\sum_{k=0}^{\tau(x)-1}(\varphi_t(x)+c)
=-t\log|dF(x)|+c\tau(x).
$$
\begin{theorem}\label{int}
Assume that $f$ admits an inducing scheme $\{S,\tau\}$ satisfying
Conditions~(H1)--(H5). Then the following statements hold:
\begin{enumerate}
\item [1.] For every $c,\ t\in\mathbb{R}$ the function $\xi_{c,t}$
satisfies Condition (P1);
\item [2.] For every $t\in\mathbb{R}$ there exists $c_t$ such that for every $c<c_t$ the potential $\xi_{c,t}$ satisfies Condition (P2) and the function $\xi_{c,t}^+$ satisfies Condition \eqref{bound1}; moreover, $P_t:=P_L(\varphi_t)$ is finite for all $t\in\mathbb{R}$;
\item [3.] There exist  $t_0=t_0(\lambda_1,\lambda_3, \gamma)<1$ and $t_1=t_1(\lambda_1,\lambda_3)>1$ such that $\xi_{c,t}$ satisfies Condition (P3) for every $t_0<t<t_1$ and every $c\in\mathbb{R}$ (the number $\gamma$ is defined in~\eqref{sng}); moreover, if
$\gamma\le\lambda_1$ then $t_0\le 0$.
\end{enumerate}
\end{theorem}
\begin{proof}
To prove the first statement we use Condition~(H5): for any
$c, t\in\mathbb{R}$, $n>0$, any cylinder $[a_0,\dots,a_{n-1}]$, and any $x,y\in J_{[a_0,\dots,a_{n-1}]}$, we have
\[
\left|\overline{\xi}_{c,t}(x)-\overline{\xi}_{c,t}(y)\right|
=|t|\left|\log\frac{|dF(y)|}{|dF(x)|}\right|
\le C\,|t|\,\lambda_2^{-n}
\]
for some constant $C>0$, thus proving the first statement.

To prove the second statement observe that
$$
\sum_{J\in S}\,\sup_{x\in J}\,\exp\,\overline{\xi}_{c,t}(x)
=\sum_{J\in S}\,e^{c\tau(J)}\,\sup_{x\in J}\,|dF(x)|^{-t}.
$$
It now follows immediately from Corollary~\ref{distbounds} that given $t\in\mathbb{R}$, there exists $c_t$ such that for every $c<c_t$ the potential $\xi_{c,t}$ satisfies Condition (P2). The finiteness of $P_t$ follows from Theorem \ref{boundedenergy} applied to the induced potential $\overline{\xi}_{c,t}$. Indeed, by Statement 1, it satisfies Condition (P1) and hence has summable variations. By Statement 2, it satisfies Condition (P2) and hence has finite Gurevich pressure. Then $P_L(\varphi_t+c)=P_t+c$ is finite and thus so is $P_t$. Now the fact that the function $\xi_{c,t}^+$ satisfies Condition \eqref{bound1} is immediate.

To establish  the remaining statements we need the following lemma.
\begin{lemma}\label{st}
We have that $P_1=0$ and
$$
P_t\ge
\begin{cases}
(1-t)\log\lambda_1 &\text{ for }t\le 1;\\
(1-t)\log\lambda_3 &\text{ for }t\ge 1.
\end{cases}
$$
\end{lemma}
\begin{proof} By the Margulis-Ruelle inequality, we have for any $f$-invariant measure $\mu$,
$$
h_\mu(f)\le \int_X\,\log|df|\,d\mu
$$
and hence, $P_1\le 0$. To show the opposite inequality note that by Conditions (H1) and (H2), for any cylinder $[a_0,\dots,a_{n-1}]$ we have that $F^n(J_{[a_0,\dots,a_{n-1}]})=W$. By the mean value theorem and Conditions (H4) and (H5) there exists a constant $c_7>0$ such that for any $x\in\overline{J_{a_0}}$ we have
\[
\text{Leb}(W)\ge c_7|dF^n(x)|\text{Leb}(J_{[a_0,\dots,a_{n-1}]}).
\]
It follows from Condition (H4) that
$$
\sum_{[a_1,\ldots,a_{n-1}]}J_{[a_0,\ldots,a_{n-1}]}=J_{a_0}.
$$
Any cylinder $[a_1,\ldots,a_{n-1}]$ contains a unique fixed point
$\omega=\omega_{[a_1,\ldots,a_{n-1}]}\in[a_1,\ldots,a_{n-1}]$ and its image $x=h(\omega_{[a_0,\ldots,a_{n-1}]})$ lies in $W$ and is a periodic point for the induced map $F$. Since $\overline\varphi_1=-\log|dF|$, we have
$$
\begin{aligned}
P_G(\overline\varphi_1)&
=\lim_{n\to\infty}\frac{1}{n}\log\sum_{F^n(x)=x\in\overline{J_{a_0}}}|dF^n(x)|^{-1}\\
&\ge\lim_{n\to\infty}\frac{1}{n}\log\sum_{[a_1,\ldots, a_{n-1}]}
\frac{c_7\,\text{Leb}(\bar{J}_{[a_0,\dots,a_{n-1}]})}
{\text{Leb}(W)}\\
&\ge\lim_{n\to\infty}\frac{1}{n}\log \frac{c_7\,\text{Leb}(\overline{J_{a_0}})}{\text{Leb}(W)}=0.
\end{aligned}
$$
By Proposition~\ref{sarig}, given $\varepsilon>0$, there exists $\nu\in\mathcal{M}(F,W)$ with $\int_W\overline\varphi_1\,d\mu>-\infty$ such that
\[
h_\nu(F)- \int_W\,\log|dF|\,d\nu\ge P_G(\overline\varphi_1)-\varepsilon\ge -\varepsilon.
\]
Since $P_1\le 0$ and $\int_W\overline\varphi_1\,d\mu>-\infty$, we also have that $h_\nu(F)<\infty$ which also yields $\int_W\overline\varphi_1\,d\mu<\infty$. In view of Corollary~\ref{distbounds}, this implies $Q_\nu<\infty$, hence, ${\mathcal L}(\nu)\in\mathcal{M}_L(f,X)$. By Theorem~\ref{Kac},
\[
\begin{aligned}
P_1\ge h_{{\mathcal L}(\nu)}(f)-&\int_X\,\log|df|\,d{\mathcal L}(\nu)\\
&=\frac{h_\nu(F)-\int_W\,\log|dF|\,d\nu}{Q_{\nu}}\ge-\frac{\varepsilon}{Q_{\nu}}\ge-\varepsilon.
\end{aligned}
\]
As $\varepsilon$ is arbitrary, $P_1\ge 0$ and we conclude that $P_1=0$. Now observe that
\[
\begin{aligned}
P_t&=\sup_{\mu\in\mathcal{M}_L(f,X)}(h_\mu-t\int_X\,\log|df|\,d\mu)\\
&\ge h_{\mu_1}-t\int_X\,\log|df|\,d\mu_1
=(1-t)\int_X\,\log|df|\,d\mu_1
\end{aligned}
\]
and the desired result follows from Theorem~\ref{integrableN}.
\end{proof}
To prove the third statement of Theorem~\ref{int} observe that
\[
\begin{aligned}
\sum_{\stackrel{J\in S}{\tau(J)\ge\tau_0}}\,&\tau(J)\sup_{x\in J}\,\exp(
\xi_{c,t}^+(x)+\varepsilon_0 \tau(x))\\
&=\sum_{\stackrel{J\in S}{\tau(J)\ge\tau_0}}\,\tau(J)e^{(-P_t+\varepsilon_0)\tau(J)}
\sup_{x\in J}|dF(x)|^{-t}=:T_t
\end{aligned}
\]
Set
$$
t_1:=\log\lambda_3(\log\frac{\lambda_3}{\lambda_1})^{-1}>1.
$$
To prove the finiteness of $T_t$ consider the following three cases:

{\bf Case I:} $1\le t<t_1$. Then $-t\log\lambda_1-P_t<0$ and
Condition~(H4) and Corollary~\ref{distbounds} yield
\[
\begin{aligned}
T_t &\le (c_3)^t\sum_{n\ge\tau_0}\,ne^{(-P_t+\varepsilon_0)n}\sum_{\tau(J)=n}\,|J|^{t-1}|J|\\
&\le (c_3)^t\sum_{n\ge\tau_0}\,n(e^{-P_t+\varepsilon_0}
\lambda_1^{-t})^n<\infty
\end{aligned}
\]
for any $0\le\varepsilon_0<t\log\lambda_1+P_t$.

{\bf Case II:} $0\le t\le 1$. Jensen's inequality yields
\[
\begin{aligned}
T_t&\le (c_3)^t\sum_{n\ge\tau_0}\,ne^{(-P_t+\varepsilon_0)n} S(n)^{1-t}
\bigl(\sum_{\stackrel{J\in S}{\tau(J)=n}}|J|\bigr)^t\\
&\le c_6^{1-t}(c_3)^t\sum_{n\ge\tau_0}\,
n(e^{(-P_t+\varepsilon_0)}\gamma^{1-t}\lambda_1^{-t})^n<\infty
\end{aligned}
\]
for any $0\le\varepsilon_0<(t-1)\log\gamma+t\log\lambda_1+P_t$.
By Lemma~\ref{st} the right hand side is positive for all
\begin{equation}\label{threestars}
t>1-\frac{\log\lambda_1}{\log\gamma}.
\end{equation}
This proves the statement for
$1-\frac{\log\lambda_1}{\log\gamma}< t\le 1$. If $\gamma\ge\lambda_1$, set
$0\le t_0:=1-\frac{\log\lambda_1}{\log\gamma}<1$. Otherwise,
$1-\frac{\log\lambda_1}{\log\gamma}$ is negative so Condition (P3) is satisfied for all values of $0\le t\le 1$. In this case $t_0=0$.

{\bf Case III:} $t\le 0$. Then
\[
\begin{aligned}
T_t&\le c_3^t \sum_{n\ge\tau_0}\,ne^{(-P_t+\varepsilon_0)n}
S(n)\lambda_3^{-tn}\\
&\le c_3^t c_6 \sum_{n\ge\tau_0}\,n(e^{(-P_t+\varepsilon_0)}
\gamma\lambda_3^{-t})^n<\infty
\end{aligned}
\]
for any $0\le\varepsilon_0<-\log\gamma+t\log\lambda_3+P_t$. Again, by Lemma~\ref{st}, the right hand side is positive provided
\[
t\ge
\log\frac{\gamma}{\lambda_1}
(\log\frac{\lambda_3}{\lambda_1})^{-1}=:t_0
\]
and $t_0<0$ if $\gamma<\lambda_1$. This completes the proof of the third statement.

\end{proof}

We now establish existence and uniqueness of equilibrium measures.

\begin{theorem}\label{tail}
Let $f$ be a $C^1$ map of a compact interval admitting an inducing scheme $\{S,\tau\}$ satisfying Conditions~(H1)--(H5). There exist constants $t_0$ and $t_1$ with $t_0<1<t_1$ such that for every $t_0<t<t_1$ one can find a measure $\mu_t\in\mathcal{M}_L(f,X)$ satisfying:
\begin{enumerate}
\item [1.] $\mu_t$ is the unique equilibrium measure (with respect to the class of liftable measures $\mathcal{M}_L(f,X)$) for the function
$\varphi_t=-t\log |d f|$;
\item [2.]  $\mu_t$ is ergodic, has exponential decay of correlations and satisfies the CLT for the class of functions whose induced functions are bounded H\"older continuous;
\item [3.]  assume that the inducing scheme $\{S,\tau\}$ is such that
$\gamma<\lambda_1$, then $t_0\le 0$ and $\mu_0$ is the unique measure of maximal entropy (with respect to the class of liftable measures $\mathcal{M}_L(f,X)$).
\end{enumerate}
\end{theorem}
\begin{proof}
Statements 1 and 3 follow directly from Theorems~\ref{remark} and \ref{int}. For Statement 2 we only  need to prove that the potential $\psi_t:=\overline{\varphi_t-P_t}$ has exponential tail with respect to the measure $i(\mu_t)=\nu_{\psi_t}$ (see Condition (P4)). By Theorem~\ref{int}, $\psi_t=\xi^+_{c, t}$ satisfies Condition (P3)
for every $t_0<t<t_1$. As $i(\mu_t)$ is a Gibbs measure there exist constants $c_8>0,\ K >0$, and $0<\theta<1$ such that
$$
\sum_{\tau(J)\ge n}\nu_{\psi_t}(J)\le c_8\sum_{\tau(J)\ge n}
\exp(\sup_{x\in J}(\overline\varphi_t(x)-P_t \tau(x)))\le K\theta^n.
$$
The statement now follows from Theorem~\ref{ergodic}.
\end{proof}

We conclude this section with the following statement.
\begin{theorem}\label{zeroentropy}
Let $f$ be a $C^1$ map of a compact interval admitting
an inducing scheme $\{S,\tau\}$ satisfying Conditions~(H1)--(H5).
Assume there exists $c_9>0$ such that for every
$\mu\in\mathcal{M}(f,I)$ with $h_\mu(f)=0$ the Lyapunov exponent
$\lambda(\mu)> c_9$. Then there exist $a>0$ and $b>0$ such that measures $\mu\in\mathcal{M}(f,I)$ with $h_\mu(f)=0$ cannot be equilibrium measures for the potential function $\varphi_t$ with
$-a<t<1+b$.
\end{theorem}
\begin{proof}
Assuming the contrary let $\mu\in\mathcal{M}(f,I)$ with
$h_\mu(f)=0$ be an equilibrium measure for $\varphi_t$. For $t>0$
\[
P_t\le h_\mu(f)-t\int_X\log|df(x)|\,d\mu(x)=-t\lambda(\mu)<-tc_9.
\]
On the other hand, since $P_t$ is decreasing we have $P_t\ge P_1=0$ for $0\le t\le 1$ leading to a contradiction. By continuity, there exists $b>0$ such that the statement also holds for $1\le t<1+b$. Since $I$ is compact, the Lyapunov exponent of a $C^1$ map $f$ is bounded from above and the same reasoning leads to a contradiction for $t>-a$ for some positive~$a$.
\end{proof}


\section{Unimodal Maps}\label{examples}


When looking for examples illustrating our theory, we may choose to
stress two different points of view: on the one hand one can strive for
the largest possible set of functions, which admit a unique equilibrium
measure; on the other hand, one might be interested in obtaining as
many potentials as possible. For unimodal maps we will give examples
in both directions.

\subsection{Definition of unimodal maps} Let $f:[b_1, b_2]\to [b_1, b_2]$ be a $C^3$ interval map with exactly one non-flat critical point (without loss of generality assumed to be $0$). Suppose $f(x)=\pm|\theta(x)|^l+f(0)$ for some local $C^3$ diffeomorphism $\theta$ and some $1<l<\infty$ (the order of the critical point). Such a map $f$ is called {\it unimodal} if $0\in (b_1, b_2)$, the derivative $df/dx$ changes signs at $0$ and
$f(b_1), f(b_2)\in\{b_1, b_2\}$.  An $S$-unimodal map is a unimodal map with negative Schwarzian derivative (for details see, for instance \cite{dMvS1}).
\begin{remark}
The negative Schwarzian derivative assumption is not necessary to prove distortion bounds for $C^3$-unimodal maps with no neutral periodic cycles \cite{Koz1} (and even $C^{2+\eta}$ unimodal maps, see \cite{Todd07}), or for $C^3$ multimodal maps \cite{VargasvanStrien04}. However, the negative Schwarzian derivative assumption avoids the simultaneous occurrence of various types of attractors in the unimodal case, so for the sake of clarity we rather assume it than restrict the statements of our theorems to the basins of the attractors.
\end{remark}

For any $x\in [b_1,b_2],\, x\neq 0$ there exists a unique point denoted by $-x\neq x$ with $f(x)=f(-x)$. If $f$ is symmetrical with respect to $0$, the minus symbol corresponds to the minus sign in the usual sense. Note that $-b_1=b_2$ so without loss of generality, we may assume that the fixed boundary point is $b:=b_2>0$ and
$f:I:=[-b,b]\to I$. If there are no non-repelling periodic cycles there exists another fixed point $\alpha$ with $f'(\alpha)<-1$ and
$0\in(\alpha, b)$. Let $\alpha^1$ denote the
unique point in $(-b, \alpha)$ for which $f(\alpha^1)=-\alpha$ and let $A=(\alpha,-\alpha)\subseteq (\alpha^1,-\alpha^1)=\hat{A}$.  An open interval $J$ is called \emph{regular} of order $\tau(J)\in\mathbb{N}$ if $f^{\tau(J)}(J)=A$ and there exists an open interval
$\hat{J}\supseteq J$ such that the map $f^{\tau(J)}|{\hat{J}}:\hat{J}\to\hat{A}$ is a diffeomorphism onto $\hat{A}$. A regular interval $J$ is called {\it maximal regular} if for every regular interval $J'$ with $J'\cap J\ne\emptyset$ we have
$J'\subseteq J$. Any two maximal regular intervals are disjoint but their closures may intersect at a boundary point. Denote by $Q$ the collection of maximal regular intervals, which are strictly contained in~$A$, and set
$$
\mathcal{\hat W}:=\bigcup_{J\in Q}\,J \,\text{ and }\,
W:=\bigcap_{n\ge 0}\, F^{-n}(A),
$$
where $F:\mathcal{\hat W}\to A$ is the induced map given by $F(x)=f^{\tau(x)}$, $x\in\mathcal{\hat W}$. Note that $W$ is the maximal $F$-invariant subset in $A$, i.e., $F^{-1}(W)=W$ and that $W\subset\hat{\mathcal{W}}$. We define
$$
S:=\{J\cap W:J\in Q\}, \quad \tau(J\cap W)=\tau(J).
$$
\subsection{Strongly regular parameters and the Collet-Eckmann condition} We consider a one-parameter family of unimodal maps $\{f_a\}$, which depends smoothly on the parameter $a$. Let
\begin{equation}\label{n0}
N_0=N_0(a):=\min\{n\in\mathbb{N} : |f_a^n(0)|<|\alpha|\}
\end{equation}
and let $F_a(0):=f_a^{N_0}(0)$. Define
$N_k:=N_{k-1}+\tau(F^{k}_a(0))$ for $k\ge 1$, where
$F_a^{k-1}(0):=f_a^{N_{k-1}}(0)$ (provided that
$F_a^{k-1}(0)\in{\mathcal{\hat W}}$). We call a parameter $a$ \emph{strongly regular} if for all $k\in\mathbb{N}$ we have
\begin{equation}\label{strongregular}
F^k_a(0)\in\mathcal{\hat W}\qquad \mbox{ and }\qquad
\sum_{\stackrel{1\leq i\leq k,}{\tau(F_a^i(0))\ge\overline{M}}}
\tau(F_a^i(0))< \rho k,
\end{equation}
where $\overline{M}=\overline{M}(N_0)$ and $\rho=\rho(N_0)$ are constants satisfying
$$
\log^2 N_0<\overline{M}<\frac{2}{3}N_0 \text{ and }
\overline{M} 2^{-\overline{M}}\ll\rho\ll 1.
$$
We denote by $\mathcal{A}$ the set of all strongly regular parameters. Observe that for any $a\in\mathcal{A}$, the first return time of the critical point to the interval $A=(-|\alpha|, |\alpha|)$ is $N_0$. Given an integer $N>0$, we denote by
$$
\mathcal{A}(N)=\{a\in\mathcal{A}: N_0(a)=N\}.
$$
Note that $\mathcal{A}=\bigcup_{N>0}\mathcal{A}(N)$.

Recall that a unimodal map satisfies the Collet-Eckmann condition if there exist constants $c>0$ and $\vartheta>1$ such that for every $n\ge 0$,
$$
|Df^n(f(0))|>c\ \vartheta^{n}.
$$
It is shown in Corollary 5.5 of \cite{Senti3} that a unimodal map $f_a$ with $a\in\mathcal{A}$ satisfies the Collet-Eckmann condition.

\subsection{Inducing schemes for unimodal maps} From now on we assume that $\{f_a\}$ is a one-parameter family of unimodal maps  with non-flat critical point in a neighborhood of a pre-periodic parameter $a^*$, that is there exists $L\in\mathbb{N}$ such that $f^L_{a^*}(0)=:x^*$ is a non-stable periodic point of period $p$. The (periodic) point $\chi(a)=f^p_a(\chi(a))$ of period $p$ for the map $f_a$ such that $\chi(a^*)=f^L_{a^*}(0)=x^*$ is called the \emph{continuation} of the point $x^*$.
Following \cite{TTY2} we call such a family of unimodal maps \emph{transverse} provided
$$
\frac{d}{da}f^L_{a^*}(0)\neq\frac{d}{da}\chi(a^*).
$$
\begin{theorem}\label{regproposition}
Let $\{f_a\}$ be a transverse one-parameter family of unimodal maps at a pre-periodic parameter $a^*$ and $\mathcal{A}$ the set of strongly regular parameters. Then
\begin{enumerate}
\item $a^*$ is a Lebesgue density point of $\mathcal{A}$, i.e.,
$$
\lim_{\varepsilon\to 0}\frac{\text{Leb}([a^*, a^*-\varepsilon]\cap\mathcal{A})}{\varepsilon}=1;
$$
moreover, there exists $T>0$ such that
$\text{Leb}(\mathcal{A}(N))>0$ for all $N\ge T$;
\item for any $f_a$ with $a\in\mathcal{A}$ the pair $\{S,\tau\}$ forms an inducing scheme satisfying Conditions (H1)--(H5). Moreover, $\text{Leb}(A\setminus W)=0$ where $W=W(a)$ is the base.
\end{enumerate}
\end{theorem}
\begin{proof}
The set of strongly regular parameters has a Lebesgue density point at $a=-2$ for the quadratic family \cite{JCY3} (see also Propositions 4.2.1 and 4.2.15 of \cite{Senti2}). A simple modification of the arguments presented there allows one to prove the same result for a transverse one-parameter family of unimodal maps at any pre-periodic parameter. The first statement follows.

Condition~(H1) follows from the definition of the collection $S$ of basic elements and Condition (H2) holds, since the induced map $F$ is expanding. To prove Condition (H3) consider a point $\omega\in S^{\mathbb N}\setminus h^{-1}(W)$. There exists $n$ such that the point $h(\sigma^n(\omega))$ is one of the end points of a maximal regular interval. It follows that the set
$S^{\mathbb N}\setminus h^{-1}(W)$ is at most countable and hence cannot support a measure, which is positive on open sets. Condition~(H4) is proven in \cite{JCY3} and Proposition 6.3 of \cite{Senti3} ( see also \cite{Senti2}) for the quadratic map. It is also shown there that the base $W$ has full Lebesgue measure in $A$. Similar arguments work for any transverse family of unimodal maps using the fact that any
non-renormalizable map of a full unimodal family is
quasi-symmetrically conjugated to a map in the quadratic family (see \cite{JakSw2}). Condition (H5) follows from Koebe's Distortion Lemma (see for example, \cite{dMvS1}).
\end{proof}
We now show that the inducing scheme $\{S,\tau\}$ satisfies Condition (H6), i.e., the number $S(n)$ of elements $J\in S$ with inducing time $\tau(J)=n$ grows subexponentially with $n$. By \cite[Proposition~2.2]{Senti3}, the partition elements of $\mathcal R$ (see Condition (H2)) of higher order are pre-images of partition elements of lower order. Hence in order to control $S(n)$, we need to control the number of intervals of lower order, which give rise to intervals of higher order. To do this we need to introduce some extra notation following \cite{Senti3}.

Denote by $J(k)$ the maximal regular interval containing $F_a^k(0)$ and by $B(k)$ the regular interval containing $F_a(0)$ for which
$f_a^{N_{k-1}-N_0}(B(k))=J(k-1)$. Let $A(k)$ be the largest interval around $0$ for which $f_a^{N_0}(A(k))\subseteq B(k)$ and let $L(k)$ be the largest regular interval in $\hat{B}(k)~\setminus~B(k)$ for which
$f_a^{N_0}(\partial A(k))$ is a boundary point. Also denote by $\hat{A}(k)$ the largest interval containing $0$ for which
$f_a^{N_0}(\hat{A}(k))\subseteq \overline{B(k)\cup L(k)}$. Finally, let $\xi_{k-1}:=f_a^{N_{k-1}}(\partial\hat{A}(k))$ and
$$
\mathcal{K}_k:=\{\text{regular intervals }J:\ f_a^{N_k}(0)\in \hat{J} \text{ and }
J\not\subseteq [\xi_k, \beta]\}.
$$
By Proposition~3.1 of \cite{Senti3}, pre-images $F_a^{-k}(J)$ of elements $J\in S$ are also elements of $S$, unless either $F^k_a(0)\in J$ or
$F^k_a(0)\in \hat{J}\setminus J$. In the first case, $J=J(k)$ and in the
second case, $J\in\mathcal{K}_k$. Since $f_a^{N_k}|\hat{A}(k)\setminus\text{int}(A(k+1))$ has two monotone
branches, for any element $J\in S$ and any $k\in\mathbb{N}$ the set
$\hat{A}(k)\setminus\text{int}(A(k+1))$ contains at most two intervals
in $S$ (of order $\tau(J)+N_k$) whose image under $f_a^{N_k}$ is
$J$. Also, for each $J'\in\mathcal{K}_k$  there are at most two intervals in $S$ (of order $\tau(J)+\tau(J')+N_k$) whose image under $f_a^{N_k+\tau(J')}$ is $J$. For strongly regular parameters Proposition~2.6 in \cite{Senti3} implies that for any interval $J'\in\mathcal{K}_k$ we have $1\le\tau(J')<\overline{M}$ if $k\le\lfloor\rho^{-1}\overline{M}\rfloor$ and $1\le\tau(J')<\rho k$ otherwise (the brackets $\lfloor\cdot\rfloor$ denote the integer part). Since all intervals in $\mathcal{K}_k$ have different order, we have that $\text{Card}\,(\mathcal{K}_k)\le\max\{1,\rho k\}$.

\begin{theorem}\label{number}
For any $\gamma>1$ there exists $c=c_\gamma>0$ and an integer $N_0>0$ such that for any $a\in\mathcal{A}(N_0)$ we have
$$
S(n)<c_\gamma\gamma^n.
$$
\end{theorem}
\begin{proof}
Observe that $S(n)=0$ for $n\in\{0,1,N_0-1, N_0\}$ and $S(n)\le
2$ for $2\le n\le N_0-2$ (see \cite[Proposition 2.2]{Senti3}). Note that
$N_0-1+2i\le N_i$ and
$$
2\gamma^{-N_0+1}\sum_{i=0}^\infty(2+\rho i)\gamma^{-2i}<1
$$
for sufficiently large $N_0$. By induction, we conclude that if
$N_k<n\le N_{k+1}$ then
\begin{multline*}
S(n)\le 2 \sum_{i=0}^k\Big(S(n-N_i)+\sum_{J'\in\mathcal{K}_i}S(n-N_i-\tau(J'))\Big)\le\\
2c_\gamma\gamma^n \sum_{i=0}^k
\Big(\gamma^{-N_i}+\sum_{J'\in\mathcal{K}_i}\gamma^{-N_i-\tau(J')}\Big)\le
2c_\gamma\gamma^{n-N_0+1}\sum_{i=0}^k(2+\rho i)\gamma^{-2i}\le\\
2c_\gamma\gamma^{n-N_0+1}\sum_{i=0}^\infty(2+\rho i)\gamma^{-2i}<c_\gamma\gamma^{n},
\end{multline*}
The desired result follows.
\end{proof}

\subsection{The liftability property for unimodal maps}\label{liftab}
We establish liftability of measures $\mu\in\mathcal{M}(f,X)$ of positive entropy, which give positive weight to the base $W$.
For a multidimensional extension of this Theorem see \cite{PSZ07}. We fix a map $f=f_a$ where $a$ is a strongly regular parameter.
\begin{theorem}\label{Bruin}
Assume that $\mu\in\mathcal{M}(f,X)$ and $h_\mu(f)>0$.
Then there exists $\nu\in\mathcal{M}(F,W)$ with ${\mathcal L}(\nu)=\mu$, i.e., $\mu\in\mathcal{M}_L(f,X)$.
\end{theorem}
\begin{proof}
Consider the Markov extension $(\underline I,\underline f)$ (also
called the Hofbauer-Keller tower) of the map $f$ (see \cite{Hof1}).
Define
$$
\underline{F}|_{\underline{\pi}^{-1}(J)}:= \underline{f}^{\tau(J)}|_{\underline{\pi}^{-1}(J)}, \quad J\in S
$$
and then
$$
\underline{A}:=\displaystyle{\bigcup_{k\ge 1} \underline{F}^k(\text{inc}\,(\bigcup_{J\in Q}J))},
$$
where $\text{inc}$ denotes the inclusion of the interval into the
first level of $\underline{I}$ and $\underline{\pi}$ the projection
from $\underline{I}$ onto the interval $I$.
By \cite{Keller3}, any $f$-invariant measure $\mu$ with $h_\mu(f)>0$
can be lifted to a measure $\underline{\mu}=\underline{\pi}^*\mu$ on
the Markov extension.

By \cite{Bru2}, if the inducing scheme is \emph{naturally extendible},
then the induced map $F$ is conjugated to the first return time map of
$\underline{A}$ via the projection map $\underline\pi$. It is easy to show that the inducing scheme constructed in Theorem~\ref{regproposition} is naturally extendible, since the intervals considered are maximal with respect to inclusion. Using the arguments in \cite[Theorem 6]{Bru2} (see also \cite{PSZ07}) we show that if $\mu\in\mathcal{M}(f, X)$ with $h_\mu(f)>0$ and $\underline{\mu}(\underline{A})>0$, then $\mu\in\mathcal{M}_L(f, X)$ as follows. Kac's formula for the first return time map $\underline{F}=\underline{f}^{\underline{R}}$ (where $\underline{R}$ is the first return time) of $\underline{A}$ to itself with $\nu=\underline{\nu}\circ\underline{\pi}^{-1}$ for the $\underline{F}$-invariant probability measure $\underline{\nu}$ yields
$$
\int\tau\,d\nu=\int \underline{R}\,d\underline{\nu}=
\frac{\underline{\mu}(\bigcup_{k\ge 0}\underline{f}^k(\underline{A}))}
{\underline{\mu}(\underline{A})}<\infty.
$$
Note that
$$
\underline{\mu}=\lim_{n\to\infty}\frac1n\sum_{k=0}^{n-1}\,\underline{\mu}_1\circ\underline{f}^k,
$$
where $\underline{\mu}_1\circ\underline{\pi}^{-1}=\mu$, and
we obtain $\nu\ll\mu$. By Zweim\"uller's dichotomy rule
\cite[Lemma2.1]{Zweimuller7}, we obtain that ${\mathcal L}(\nu)=\kappa\cdot\mu$ for some $\kappa>0$. Normalizing $\nu$ if necessary one has that $\mu\in\mathcal{M}_L(f, X)$.

To prove that $\underline{\mu}(\underline{A})>0$ for any
$\mu\in\mathcal{M}(f, X)$ it suffices to show that
\begin{equation}\label{liftabilitycondition}
\underline{\pi}^{-1}(X)\subseteq \bigcup_{k\ge
0}\underline{f}^{-k}(\underline{A})\pmod{\underline{\mu}}.
\end{equation}
Indeed, in view of
\eqref{liftabilitycondition}, the assumption that
$\underline{\mu}(\underline{A})=0$ leads to the following
contradiction:
$$
1=\mu(X)=\underline{\mu}\circ\underline{\pi}^{-1}(X)
\le\sum_{k\ge 0}\,\underline{\mu}(\underline{f}^{-k}
(\underline{A}))=\sum_{k\ge 0}\,\underline{\mu}(\underline{A})=0.
$$
In order to establish \eqref{liftabilitycondition} for the inducing scheme constructed in Theorem~\ref{regproposition} observe that by (H2), any point $x\in X$ has a basis of neighborhoods, which are sent
diffeomorphically by some iterates of $f$ onto $\hat{A}$ (i.e., the
extension of $A$). Denote the (countable) set of boundary points of
$\underline{I}$ by $\partial\underline{I}$. Without loss of generality we may assume that $\mu$ has no atoms and thus $\underline{\mu}(\partial\underline{I})=0$. By the Markov property of $(\underline{I}, \underline{f})$, any point
$\underline{x}\in\underline{\pi}^{-1}(X)\setminus\partial\underline{I}$
has a basis of neighborhoods $\underline{U}\subset \underline{\hat{U}}$ such that for some integer $k$ and some level $\underline{D}_\ell$ of $\underline{I}$ we have
$$
\underline{\pi}\circ\underline{f}^{k}(\underline{U})
=\ A\ \subset\, \underline{\pi}\circ\underline{f}^{k}(\underline{\hat{U}})
= \hat{A}\,\subseteq\ \underline\pi\,(\underline{D}_\ell)
$$
(recall that the $\ell$-th level of $\underline{I}$ is the image under $\underline{f}^\ell$ of the maximal interval of monotonicity of $f^\ell$).
Therefore we are left to show that for any
$\underline{\hat{A}}_{\ell}\in\underline{\pi}^{-1} (\hat{A})\cap\underline{D}_\ell$ we have
$$
\exists\ \underline{A}_{\ell}\subset \underline{\hat{A}}_\ell,\ \underline{A}_{\ell}\in\underline{A} \quad\iff\quad\hat{A}\subseteq\ \underline{\pi}\,(\underline{D}_\ell).
$$
For our partition the {\em ``$\ \Rightarrow$''} direction follows from the arguments of \cite[Lemma2]{Bru2}.

We are left to prove that if
$\hat{A}\subseteq\,\underline{\pi}\,(\underline{D}_\ell)$ then there
exists some set $J\in S$ and some integer $k$ such that
$\underline{F}^k\,(\text{inc}\, (J))=\underline{A}_\ell$. Recall
that $\underline{D}_\ell=[c_{\ell-i}, c_\ell]$, where $c_n:=f^n(0)$
and
$i:=\max_{1\le j<\ell}\{ 0\in\underline{D}_j\}<\ell$. Denote by
$c_{-n}$ the $n$-th pre-image of the critical point, which lies closest
to the critical point. We have that for $0\le k\le l$
$$
0\not\in f^k(]c_{-\ell}, 0[),\quad
\underline{f}^\ell(\text{inc}\,( [c_{-(\ell-i)},0]))=\underline{D}_\ell
$$
and there exist $J_\ell\subset\hat{J}_\ell\subset [c_{-(\ell-i)}, 0]$ for which
$$
\underline{f}^\ell|_{\text{inc}\,(J_\ell)}(\text{inc}\,(J_{\ell}))
=\underline{A}_{\ell}\subset
\underline{f}^\ell|_{\text{inc}\,(\hat{J}_\ell)}(\text{inc}\,(\hat{J}_{\ell}))
=\underline{\hat{A}}_{\ell}
\subseteq\underline{D}_\ell.
$$
If $J_\ell\in S$ then
$\underline{F}(\text{inc}\,(J_\ell))=\underline{A}_\ell$ and
$\underline{A}_\ell\in\underline{A}$. Otherwise, $J_\ell\subset J\in S$. Again, there are no pre-images of the critical point of order less than $\tau(J)$ in $J$ or between $J$ and the critical point, so
$\underline{\pi}(\underline{f}^{\tau(J)}(\text{inc}\,J))=A$ and $\underline{f}^{\tau(J)}(\text{inc}\,J)\in\underline{A}$. Again, if
$f^{\tau(J)}(J_\ell)\in S$ then
$$
\underline{F}^2(\text{inc}\,(J_\ell))
=\underline{f}^{\tau(J)+\tau(J_\ell)}(\text{inc}\,(J_\ell))
=\underline{A}_\ell
$$
and $\underline{A}_\ell\in\underline{A}$. Inductively, this shows that there exists $J$ for which
$\underline{F}^{k}(\text{inc}\,(J))=\underline{A}_\ell$. Hence,
$\underline{A}_\ell\in\underline{A}$. This completes the proof.
\end{proof}
\begin{remark} If the inducing scheme constructed in Theorem~\ref{regproposition} is refined according to Remark~\ref{refinement}, the new inducing scheme $\{S',\tau'\}$ is no longer naturally extendible. However, one can express $\{S',\tau'\}$ as an inducing scheme over $(W, F)$. Namely, for each element $J'\in S'$ with $J'\subset J\in S$, set $\tau'(J'):=n(J)+1$ where $n(J)\ge 0$ is the number of times $J$ needs to be refined to obtain $J'$. We then have $F'(x):=F^{\tau'(J')}(x)$ for $x\in J'$. Since the refinement of Remark~\ref{refinement} is finite there exists a uniform bound on all $n(J)$ and
$$
\sum n(J')\nu(J')<\infty,
$$
by \cite[Theorem 1.1]{Zweimuller7}, hence $\nu\in\mathcal{M}_L(W, F)$. In other words, there exists an $F'$-invariant probability measure $\nu'$ on $W'$ such that ${\mathcal L}(\nu')=\nu$ and therefore,
$$
{\mathcal L}({\mathcal L}(\nu'))={\mathcal L}(\nu)=\mu\in\mathcal{M}_L(X, f).
$$
\end{remark}

We now prove that for strongly regular parameters equilibrium measures must give positive weight to the base $W$.

\begin{theorem}\label{nunonzero}
Let $\{f_a\}$ be a transverse one-parameter family of $S$-unimodal maps with non-flat critical point in a neighborhood of a pre-periodic parameter $a^*$. There exists $N_0$ such that for every $n\ge N_0$ and every $a\in\mathcal{A}(n)$ there exist $t'_0=t'_0(a)<0$ and $1<t'_1=t'_1(a)$ such that for any $t'_0<t<t'_1$ we have that
$$
\sup_{\stackrel{\nu\in\mathcal{M}(f_a, I_a)}{\nu(W)=0}}\{h_\nu(f_a)-t\lambda(\nu)\}<
\sup_{\mu\in\mathcal{M}_L(f_a, X_a)}\{h_\mu(f_a)-t\lambda(\mu)\}.
$$
where $\lambda(\nu)=\lambda_a(\nu)=\int_I\log|df_a(x)|\,d\nu$.
\end{theorem}
\begin{proof}
In the particular case of the quadratic family, we have
$$
\dim_H(A\setminus W)=
\dim_H(\bigcup_{k\ge 0}F^{-k}(A\setminus \bigcup_{J\in S}J))<c\frac{\log N_0}{N_0}
$$
for all $a\in\mathcal{A}(N)$ (see \cite{Senti2}, \cite{Senti3}) and some constant $c\in\mathbb{R}$. By definition of $X$, any $f$-invariant Borel measure $\nu$ with $\nu(W)=0$ must satisfy $\nu(X)=0$, and by construction, if $f^k(x)\in A$ for $x\in X$ then $f^k(x)\in W$. So $X$ is the disjoint union of $W$ and of its pre-images along H\"older continuous inverse branches of $f$ (they are bounded away from the critical value) and hence
\begin{equation}\label{sentidimension}
\dim_H\nu\le\dim_H((f(0), f^2(0))\setminus X)=c\dim_H(A\setminus W)
\end{equation}
for some constant $c\in\mathbb{R}$, since the support of any $f$-invariant measure $\nu\neq\delta_\beta$ is contained in $(f(0), f^2(0))$. In particular, the Hausdorff dimension can thus be made arbitrarily small by choosing the number $N_0$ to be sufficiently large. In the general case, by \cite{JakSw2}, $f_a$ is H\"older conjugated to a quadratic map, so the Hausdorff dimension of $\nu$ can also be made arbitrarily small provided $N_0$ is sufficiently large.

We now proceed with the proof of the theorem and we argue by contradiction assuming the statement is false. Then, for every
$\varepsilon>0$ there exists an invariant Borel measure $\nu$ with
$\nu(W)=0$ and such that
$$
h_\nu(f_a)-t\lambda(\nu)\ge P_{t,a}-\varepsilon,
$$
where $P_{t,a}=P_L(\varphi_{t,a})$ is defined by \eqref{supremum}.
We first consider the case when $t<1$. Then one can choose
$0<\varepsilon<\min\{(1-t)\log\lambda_1,\log\lambda_1 \}$ and
$0<\delta<\frac{\log\lambda_1-\varepsilon}{\log\lambda_3}$ where
$\lambda_1=\lambda_1(a)$ is the constant from Condition~(H4) and
$\lambda_3=\lambda_3(a)$ is such that
$\lambda(\nu)\le\log\lambda_3$ for every $f_a$-invariant measure
$\nu$ (such a constant exists, since $f$ is $C^1$ on a compact set). Young's formula for the dimension of the measure (see \cite{LSY3}) and Lemma~\ref{st} yield for the $\epsilon$ and $\delta$ above
\begin{multline*}
\dim_H\nu=\frac{h_\nu(f_a)}{\lambda(\nu)}\ge t+\frac{P_{t,a}-\varepsilon}{\lambda(\nu)}
\ge t+\frac{(1-t)\log\lambda_1-\varepsilon}{\log\lambda_3}\ge\delta>0
\end{multline*}
for every $t$ satisfying
$$
t'_0:=\Bigl(\delta-\frac{\log\lambda_1-\varepsilon}{\log\lambda_3}\Bigr)
\Bigl(1-\frac{\log\lambda_1}{\log\lambda_3}\Bigr)^{-1} \le t\le 1.
$$
Note that $t'_0$ is negative.

We now consider the case when $t\ge 1$. Recall that for any Collet-Eckmann parameter all probability measures have a strictly positive Lyapunov exponent $\lambda(\nu)\ge \lambda_{inf, a}>0$. Choose $0<\varepsilon<\lambda_{inf}$ and $0<\delta<1-\frac{\varepsilon}{\lambda_{inf}}$. By Lemma~\ref{st}, we have that $0\le P_{t,a}\ge (1-t)\log\lambda_3$ and hence
$$
\dim_H\nu\ge t(1-\frac{\log{\lambda_3}}{\lambda_{inf}})+\frac{\log{\lambda_3}-\varepsilon}{\lambda_{inf}}\ge \delta>0$$
for every $t$ satisfying
$$
1\le t\le \Bigl(\delta-\frac{\log\lambda_3-\varepsilon}{\lambda_{inf}}\Bigr)
\Bigl(1-\frac{\log\lambda_3}{\lambda_{inf}}\Bigr)^{-1}:=t'_1.
$$
Observe that $t'_1>1$. To conclude note that one can choose the set of parameters of positive Lebesgue measure such that $N_0$ is arbitrarily large and hence the dimension of $\nu$ (see \eqref{sentidimension}) is less than $\delta$. This leads to a contradiction.
\end{proof}

\subsection{Equilibrium measures for unimodal maps}
We now summarize our results on unimodal maps, observing that they extend the results of \cite{Bru-Kel} for the parameters under consideration. The proof follows from Theorems~\ref{tail}, \ref{zeroentropy}, \ref{regproposition}, \ref{number}, \ref{Bruin} and \ref{nunonzero}.

\begin{theorem}\label{mainunimodaltheorem}
Let $\{f_a\}$ be a transverse one-parameter family of $S$-unimodal
maps with non-flat critical point in a neighborhood of a pre-periodic parameter $a^*$. Then for every $\mathcal{A}(N)$ of positive measure and every $a\in\mathcal{A}(N)$
\begin{enumerate}
\item [1.] one can find numbers $t_0=t_0(a)<0$ and $t_1=t_1(a)>1$ such that for every $t_0<t<t_1$ there exists a unique equilibrium measure $\mu_{t,a}$ for the  function
$\varphi_{t,a}(x)=-t\log |df_a(x)|$, $x\in I$, i.e.,
\begin{equation}\label{phita}
\begin{aligned}
\sup\{h_\mu(f_a)-&t\int_I\log |df_a(x)|\,d\mu\}\\
&=h_{\mu_{t,a}}(f_a)-t\int_I\log |df_a(x)|\,d\mu_{t,a},
\end{aligned}
\end{equation}
where the supremum is taken over all $f_a$-invariant Borel probability
measures.
\item [2.] the measure $\mu_{t,a}$ is ergodic, has exponential decay of correlations and satisfies the CLT for the class of functions whose
induced functions are bounded H\"older continuous. In particular, there exists a unique measure $\mu_{0,a}$ of maximal entropy and a unique absolutely continuous invariant measure $\mu_{1,a}$.
\end{enumerate}
\end{theorem}

For the purpose of obtaining the largest class of functions admitting a unique equilibrium measure for $\varphi_{t, a}(x)$, we can consider the families of maps studied by Avila and Moreira in \cite{AM3, AM2}. Let us call a smooth (at least $C^3$) unimodal map \emph{hyperbolic} if it has a quadratic critical point, has a hyperbolic periodic attractor and its critical point is neither periodic, nor pre-periodic. A family of unimodal maps is called \emph{nontrivial} if the set of parameters for which the corresponding map is hyperbolic is dense. One can also consider families of maps that depend on any number of parameters. We then obtain the following result. A parameter is called regular if the corresponding unimodal map has a hyperbolic periodic attractor.

\begin{theorem}
Let $\{f_a\}$ be a nontrivial analytic family of $S$-unimodal maps. Then for almost every non-hyperbolic parameter the corresponding map $f_a$ admits a unique equilibrium measure (with respect to the class $\mathcal{M}(f_a, X)$) for the potential $\varphi_{t,a}(x)$ for all
$t_0<t<t_1$ with some $0<t_0=t_0(a)$ and $t_1=t_1(a)>1$. The same result holds for any non-regular parameters in any generic smooth ($C^k, k=2,\ldots, \infty$) family of unimodal maps.
\end{theorem}
\begin{proof}
By \cite[Theorem~A]{AM2}, \cite{AM3}, almost every non-regular parameter of a family of unimodal maps satisfying our hypothesis also satisfies the Collet-Eckmann condition. By \cite{BLvS1}, any unimodal map, satisfying the Collet-Eckmann condition, admits an inducing scheme satisfying Conditions (H1)--(H5). The result now follows from Theorem~\ref{tail}. By \cite[Proposition 3.1]{Bru-Kel}, any invariant measure has uniformly positive Lyapunov exponent. Theorems~\ref{zeroentropy} and \ref{Bruin} then imply that the equilibrium measure can be taken with respect to the class of all measures in $\mathcal{M}(X, f_a)$.
\end{proof}
Under slightly stronger regularity conditions (satisfied for instance, if
$f$ is a polynomial map) Bruin and Keller show \cite{Bru-Kel} that measures $\mu\in\mathcal{M}(I,f_a)\setminus\mathcal{M}(X,f_a)$ cannot be equilibrium measures for the potential functions
$\varphi_{t,a}(x)$ with $t$ close to $1$.

\section{More Interval Maps.}\label{multimodalsection}
\subsection{Multimodal maps}
We follow \cite{BLvS1}. Consider a $C^3$ interval or circle map $f$
with a finite set $\mathcal{C}$ of critical points and no stable or neutral periodic point. Also assume that all critical points have the same order $\ell$, i.e., for each $c\in\mathcal{C}$ there exists a diffeomorphism
$\psi:\mathbb{R}\to\mathbb{R}$ fixing $0$ such that for $x$ close to
$c$ we have
$$
f(x)=\pm|\psi(x-c)|^\ell+f(c)
$$
where $\pm$ may depend on the sign of $x-c$. Assume  (as in \cite{BLvS1}) that
\begin{equation}\label{multimodal}
\sum_{n\in\mathbb{N}}|df^n(f(c))|^{\frac{-1}{2\ell-1}}(c)<\infty
\qquad \mbox{ for each } c\in\mathcal{C}
\end{equation}
and that there exists a sequence $\{\gamma_n\}_{n\in\mathbb{N}}$, $\gamma_n\in(0,\frac12)$ satisfying
$$
\sum_{n\in\mathbb{N}}\gamma_n<\infty
$$
and for some $\beta>0$, each $c\in\mathcal{C}$ and $n\ge 1$,
\begin{equation}\label{gamma}
\left(\gamma_n^{\ell-1}|df^n(f(c))|\right)^{-\frac1\ell}\le Ce^{-\beta n}.
\end{equation}
Let $X$ be the biggest closed $f$-invariant set of positive Lebesgue measure. This set can be decomposed into finitely many invariant subsets $X_i$ on which $f$ is topologically transitive. The following result is an easy corollary of \cite[Proposition~4.1]{BLvS1}.
\begin{theorem}\label{regproposition1}
Let $f$ be a multimodal map satisfying Conditions \eqref{multimodal} and \eqref{gamma}. Then for each $i$, the map $f|X_i$ admits an inducing scheme $\{S_i,\tau_i\}$ satisfying Conditions~(H1)--(H5). The corresponding inducing domain $W_i$ lies in a small neighborhood of  a critical point and the basic elements of the inducing scheme accumulate to the critical point.
\end{theorem}
We thus obtain the following result.

\begin{theorem} Let $f$ be a multimodal map satisfying Conditions
\eqref{multimodal} and \eqref{gamma}. Then for every $X_i$ there
exist $t_0<1<t_1$ such that for every $t_0<t<t_1$ one can find a
unique equilibrium measure $\mu_{t,i}$ on $X_i$ for the function
$\varphi_t=-t\log |df|$ with respect to the class of measures
$\mathcal{M}_L(f,X_i)$. The measure $\mu_{t,i}$ is ergodic, has exponential decay of correlations and satisfies the CLT for the class of functions whose induced functions are bounded H\"older continuous. Additionally, if $f$ satisfies the Collet-Eckmann condition (for multimodal maps), then  $\mu_{t,i}$ is the unique equilibrium measure with respect to the class of measures $\mathcal{M}(f,X_i)$.
\end{theorem}

\begin{proof}
The first part is a direct corollary of Theorem~\ref{regproposition1}. To prove that the equilibrium measure is unique with respect to all invariant measures in $\mathcal{M}(f,X_i)$, we remark that
Theorem~\ref{Bruin} holds for any piecewise continuous piecewise monotone interval map provided the basic elements of the inducing scheme accumulate to the critical point (see \cite[Section~7]{PSZ07} for details and more general results). This implies that the class $\mathcal{M}_L(f,X_i)$ includes all $f$-invariant measures on $X_i$ of positive entropy (\cite{Hof1}). By \cite[Theorem~1.2]{BLvS1}, every invariant measure has Lyapunov exponent bounded away from $0$ and hence no invariant measure of zero entropy can be an equilibrium measure for the function~$\varphi_t$.
\end{proof}
\subsection{Cusp maps}
A \emph{cusp map} of a finite interval $I$ is a map $f:\bigcup_{j}I_j\to I$ of an at most countable family $\{I_j\}_j$ of disjoint open subintervals of $I$ such that
\begin{enumerate}
\item[$\circ$] $f$ is a $C^1$ diffeomorphism on each interval $I_j:=(p_j, q_j)$, extendible to the closure $\bar{I}_j$ (the extension is denoted by $f_j$);
\item[$\circ$] the limits $\lim_{\epsilon\to 0^+}Df(p_j+\epsilon)$ and $\lim_{\epsilon\to 0^+}Df(q_j-\epsilon)$ exist and are equal to either $0$ or $\pm\infty$;
\item[$\circ$] there exist constants $K_1>K_2>0$ and $C>0$, $\delta>0$ such that for every $j\in \mathbb{N}$ and every $x,x'\in\bar{I}_j$,
\[
|Df_j(x)-Df_j(x')|< C|x-x'|^\delta \text{ if } |Df_j(x)|\,,\,|Df_j(x')|\le K_1,
\]
\[
|Df^{-1}_j(x)-Df^{-1}_j(x')|< C|x-x'|^\delta \text{ if } |Df_j(x)|\,,\,|Df_j(x')|\ge K_2.
\]
\end{enumerate}

In \cite{DobbsPhD}, it is shown that certain cusp maps admit inducing schemes.

\begin{theorem}\label{cusp}
Let $f$ be a cusp map with finitely many intervals of monotonicity $I_j$. Suppose $f$ has an ergodic absolutely continuous invariant probability measure $m$ with strictly positive Lyapunov exponent. Then $f$ admits an inducing schemes $\{S, \tau\}$ which satisfies Conditions~(H1)--(H3) and (H5).
\end{theorem}
\begin{proof}
Conditions~(H1), (H2), (H5) are satisfied by the definition of the Markov maps from \cite[Theorem~1.9.10]{DobbsPhD}. To prove Condition~(H3) observe that any point of $S^\mathbb{N}\setminus h^{-1}(W)$ is eventually mapped onto an endpoint of one of the domains of the Markov map. Since these domains are intervals, the set of all endpoints is a countable set, and so the set $S^\mathbb{N}\setminus h^{-1}(W)$ cannot support a measure which is positive on open sets, proving Condition~(H3).
\end{proof}

By definition, for cusp maps one cannot expect to obtain an upper bounds on the derivatives of the induced map and of the Lyapunov exponent of liftable measures using compacity arguments as in Corollary~\ref{distbounds} and Theorem~\ref{integrableN}. However, since this upper bound is only used to extend the range of values of $t$ for which our theorems hold, one can nonetheless obtain statements on the existence of a unique equilibrium measure associated to the potential $-t\log |df|$, albeit for a smaller range of values $t$. Theorem~\ref{int} now becomes as follows.
\begin{theorem}\label{int1}
Assume that the cusp map $f$ admits an inducing scheme $\{S,\tau\}$ satisfying Conditions~(H1)--(H5). Then the following statements hold:
\begin{enumerate}
\item [1.] For every $c,\ t\in\mathbb{R}$ the function $\xi_{c,t}$ satisfies Condition (P1);
\item [2.] For every $t\ge0$ there exists $c_t$ such that for every $c<c_t$ the potential $\xi_{c,t}$ satisfies Condition (P2) and the function $\xi_{c,t}^+$ satisfies Condition \eqref{bound1}; moreover, $P_t:=P_L(\varphi_t)$ is finite for all $t\ge0$;
\item [3.] There exist $t_0^*=t_0^*(\lambda_1)<1$ and $t_1^*=t_1^*(\lambda_1)>1$ such that $\xi_{c,t}$ satisfies Condition (P3) for every $t_0^*<t<t_1^*$ and every $c\in\mathbb{R}$;
\end{enumerate}
\end{theorem}
\begin{proof}
The proof of parts 1 and 2 follow as in Theorem~\ref{int} (although Statement 2 now only holds for non-negative values of $t$). To prove Statement 3, observe that $P_1\ge0$ by \cite[Theorem~1.9.12]{DobbsPhD}, and so by continuity, there exist $t_0^*=t_0^*(\lambda_1)<1$ and $t_1^*=t_1^*(\lambda_1)>1$ such that (P3) holds for every $t_0^*<t<t_1^*$.
\end{proof}

For the inducing scheme constructed in Theorem~\ref{cusp} the liftability problem is solved in \cite[Corollary~7.5]{PSZ07}: for cusp maps every measure of positive entropy which gives positive weight to the base of the inducing scheme is liftable.

Also, one should note that while applying our results to cusp maps Condition~(H4) may not hold in general and so we must assume it. Combining this result with Theorems~\ref{tail} and \ref{int1} yield the following statement.

\begin{theorem}
Let $f$ be a cusp map with finitely many intervals of monotonicity, which admits an ergodic absolutely continuous invariant probability measure $m$ with strictly positive Lyapunov exponent. Additionally assume that Condition~(H4) is satisfied for the associated inducing schemes $\{S, \tau\}$. Then there exist $t_0<1<t_1$ such that there is a unique equilibrium measure $\mu$ (with respect to the class of all invariant measures) with $\mu(W)>0$ (where $W$ is the domain of the inducing scheme) associated to the potential function $-t\log |df|$ for all $t_0<t<t_1$. This measure is ergodic, has exponential decay of correlations and satisfies the central limit theorem for the class of functions whose induced functions are bounded H\"older continuous.
\end{theorem}

\bibliographystyle{alpha}
\bibliography{biblio1_sam}
\end{document}